\documentclass[12pt]{article}%
\usepackage{amsmath}
\usepackage{amsfonts}
\usepackage{amssymb}
\usepackage{graphicx}%
\newtheorem{theorem}{Theorem}

\newtheorem{definition}[theorem]{Definition}

\newtheorem{lemma}[theorem]{Lemma}

\newtheorem{remark}[theorem]{Remark}

\newenvironment{proof}[1][Proof]{\noindent\textbf{#1.} }{\ \rule{0.5em}{0.5em}}
\begin{document}

\title{Asymptotic analysis of the Krawtchouk polynomials by the WKB method}
\author{Diego Dominici \thanks{e-mail: dominicd@newpaltz.edu}\\Department of Mathematics\\State University of New York at New Paltz\\75 S. Manheim Blvd. Suite 9\\New Paltz, NY 12561-2443\\USA}
\maketitle

\begin{abstract}
We analyze the Krawtchouk polynomials $K_{n}(x,N,p,q)$ asymptotically. We use
singular perturbation methods to analyze them for $N\rightarrow\infty,$ with
appropriate scalings of the two variables $x$ and $n$. In particular, the WKB
method and asymptotic matching are used. We obtain asymptotic approximations
valid in the whole domain $[0,N]\times\lbrack0,N],$ involving some special
functions. We give numerical examples showing the accuracy of our formulas.

\end{abstract}

\section{Introduction}

\begin{definition}
Let $n,N\geq0$ be integers. The Krawtchouk polynomials $K_{n}(x)$ are defined
by \cite{MR51:8724},
\begin{equation}
K_{n}(x)=\sum_{k=0}^{n}\binom{x}{k}\binom{N-x}{n-k}q^{k}\left(  -p\right)
^{n-k}\label{Kndef}%
\end{equation}
where%
\begin{equation}
0<p,q<1\qquad p+q=1.\label{p,qdef}%
\end{equation}
The \textit{binary Krawtchouk polynomials} are the special case with $p=1/2=q$
\cite{MR2002d:33015}, \cite{MR2002d:33016}.
\end{definition}

The Krawtchouk polynomials \cite{MR92m:33019} are one of the families of
classical orthogonal polynomials of a discrete variable
\cite{koekoek94askeyscheme}. They satisfy the orthogonality relation
\[
\sum_{k=0}^{N}K_{i}(k)K_{j}(k)\varrho(k)\text{=}\binom{N}{j}(pq)^{j}%
\delta_{ij}\text{,\quad}i,j=0,\ldots,N
\]
with weight function%
\[
\varrho(x)=\binom{N}{x}p^{x}q^{N-x}.
\]

Writing the Krawtchouk polynomials in the extended form $K_{n}(x)=K_{n}%
(x,N,p,q)$, we have the the symmetry formula%
\begin{equation}
K_{n}(x,N,p,q)=(-1)^{n}K_{n}(N-x,N,q,p).\label{symetry}%
\end{equation}
They also satisfy the three-term recurrence
\[
(n+1)K_{n+1}(x)+pq(N-n+1)K_{n-1}(x)+\left[  p(N-n)+nq-x\right]  K_{n}(x)=0
\]
which is the main object of our analysis.

The Krawtchouk polynomials are important in the study of the Hamming scheme of
classical coding theory \cite{MR2002a:94043}, \cite{MR97c:94015},
\cite{MR57:5408a}, \cite{MR94c:94009}, \cite{MR53:5155}, \cite{MR91h:33005}.
Lloyd's theorem \cite{MR19:465b} states that if a perfect code exists in the
Hamming metric, then the Krawtchouk polynomial must have integral zeros
\cite{MR53:5150}, \cite{MR52:5187}, \cite{MR46:6918}. Not surprisingly, these
zeros have been the subject of extensive research \cite{MR91a:33011},
\cite{MR2001h:31002} \cite{MR2002m:33019}, \cite{MR2003c:33013},
\cite{MR94k:33011}, \cite{MR97i:33005}, \cite{MR2000k:33029}.

The Krawtchouk polynomials also have applications in probability theory
\cite{MR94d:60011}, queueing models \cite{saddledisc}, stochastic processes
\cite{MR2001f:60095}, quantum mechanics \cite{MR2003h:39013},
\cite{MR98j:81070} \cite{MR2002f:33014} and biology \cite{MR53:14704}.

The asymptotic behavior of the Krawtchouk polynomials as $N\rightarrow\infty$
was studied by Sharapudinov for $x\approx Np$ and $n=O\left(  N^{1/3}\right)
$ in \cite{MR90a:33013}. He derived an approximation in terms of the Hermite
polynomials (see also \cite{MR51:8724} for a similar formula). The general
case with $x,n=O(N)$ was investigated by Ismail and Simeonov
\cite{MR2000a:33015} and a uniform asymptotic expansion was derived by Li and
Wong \cite{MR2001i:33017}, both using the saddle point method.

The purpose of this paper is to take a different approach, based on the
recurrence formula that the Krawtchouk polynomials satisfy and using singular
perturbation techniques \cite{MR87g:42044}, \cite{MR2004c:39012} to analyze
it. We scale $x=yN,$ $n=zN$ and obtain asymptotic approximations to $K_{n}(x)$
for $(y,z)\in\lbrack0,1]\times\lbrack0,1].$ Our results agree and extend those
obtained in \cite{MR2000a:33015}, \cite{MR2001i:33017} and \cite{MR2081675}.

In Section 2 we review the basic properties of the Krawtchouk polynomials. In
Sections 3-10 we obtain asymptotic expansions from the recurrence formula by
using the WKB method. We must consider twelve relevant regions of the
two-dimensional state space. In Section 11 we summarize our results and
numerically compare our approximations with the exact formula.

\section{The WKB approximation}

To analyze the recurrence
\begin{equation}
(n+1)K_{n+1}(x)+pq(N-n+1)K_{n-1}(x)+\left[  p(N-n)+nq-x\right]  K_{n}(x)=0
\label{recurrence}%
\end{equation}
subject to the boundary conditions%
\begin{equation}
K_{0}(x)=1 \label{condn=0}%
\end{equation}%
\begin{equation}
K_{N+1}(x)=\binom{x}{N+1} \label{condN+1}%
\end{equation}%
\begin{equation}
K_{n}(0)=\binom{N}{n}\left(  -p\right)  ^{n} \label{condx=0}%
\end{equation}%
\begin{equation}
K_{n}(N)=\binom{N}{n}q^{n} \label{condx=N}%
\end{equation}
for large $N,$ we introduce the scaled variables $y,z$ defined by%
\begin{equation}
x=yN,\quad n=zN,\quad0<y,z<1. \label{y,z}%
\end{equation}
We define the function $F(y,z)$ and the small parameter $\varepsilon$ by%
\begin{equation}
\varepsilon=\frac{1}{N},\quad K_{n}(x)=F(\varepsilon x,\varepsilon n)=F(y,z)
\label{F(y,z)}%
\end{equation}
and observe that $K_{n\pm1}(x)=F(y,z\pm\varepsilon).$

Substituting (\ref{y,z})-(\ref{F(y,z)}) in (\ref{recurrence}) we get%
\begin{equation}
(z+\varepsilon)F(y,z+\varepsilon)+pq(1-z+\varepsilon)F(y,z-\varepsilon
)+\left[  p(1-z)+zq-y\right]  F(y,z)=0.\label{EqF}%
\end{equation}
To find $F(y,z)$ for $\varepsilon$ small, we shall use the WKB method
\cite{MR1429619}. Thus, we consider solutions which have the asymptotic form%
\begin{equation}
F(y,z)\sim\varepsilon^{\nu}\exp\left[  \varepsilon^{-1}\psi(y,z)\right]
L(y,z).\label{ansatz}%
\end{equation}
Using (\ref{ansatz}) in (\ref{EqF}), with%
\[
\varepsilon^{-1}\psi(y,z\pm\varepsilon)=\varepsilon^{-1}\psi(y,z)\pm\psi
_{z}(y,z)+\frac{1}{2}\psi_{zz}(y,z)\varepsilon+O\left(  \varepsilon
^{2}\right)  ,
\]
dividing by $\exp\left[  \varepsilon^{-1}\psi(y,z)\right]  $ and expanding in
powers of $\varepsilon$ we obtain the \textit{eikonal equation}%
\begin{equation}
zU^{2}+\left[  p-y+z\left(  q-p\right)  \right]  U+pq(1-z)=0\label{eikonal}%
\end{equation}
and the \textit{transport equation}%
\begin{equation}
\left[  zU^{2}-pq(1-z)\right]  L_{z}+\left\{  \frac{1}{2}\left[
zU^{2}+pq(1-z)\right]  \psi_{zz}+U^{2}+pq\right\}  L=0\label{transport}%
\end{equation}
with%
\begin{equation}
U(y,z)=\exp\left[  \psi_{z}(y,z)\right]  .\label{Udef}%
\end{equation}

\subsection{The functions $\psi$ and $L$}

From (\ref{Udef}) we have%
\[
\psi(y,z)=\int\ln\left[  U(y,z)\right]  dz=z\ln(U)-\int z\frac{U_{z}}%
{U}dz=z\ln(U)-\int z\frac{U_{z}}{U}\frac{dz}{dU}dU
\]
and using (\ref{eikonal}) we get%
\[
\int z\frac{U_{z}}{U}\frac{dz}{dU}dU=\int\frac{(y-p)U-pq}{(U+q)(U-p)}%
dU=\int\left[  \frac{1}{U}+\frac{y-1}{U}-\frac{y}{U+q}\right]  dU.
\]
Hence,%
\begin{equation}
\psi(y,z)=\ln\left[  U^{z-1}\left(  U-p\right)  ^{1-y}\left(  U+q\right)
^{y}\right]  +A(y) \label{psi}%
\end{equation}
where the function $A(y)$ is still unknown.

From (\ref{transport}) we have%
\[
L(y,z)=B(y)\exp\left[  -%
{\displaystyle\int}
\frac{\frac{1}{2}\left[  zU^{2}+pq(1-z)\right]  \psi_{zz}+U^{2}+pq}%
{zU^{2}-pq(1-z)}dz\right]
\]
and from (\ref{Udef}) $\psi_{zz}=U_{z}/U.$ After changing variables from $z$
to $U,$ we obtain%
\[
L(y,z)=B(y)\exp\left\{
{\displaystyle\int}
\left[  \frac{1}{2U}+\frac{1}{U+q}-\frac{\left(  p-y\right)  U+pq}{\left(
p-y\right)  U^{2}+2pqU+pq(q-y)}\right]  dU\right\}  .
\]
Thus,%
\begin{equation}
L(y,z)=B(y)\left(  U+q\right)  \sqrt{\frac{U}{\left(  p-y\right)
U^{2}+2pqU+pq(q-y)}} \label{L1}%
\end{equation}
where $B(y)$ is to be determined.

\subsection{The function $U$}

Rewriting (\ref{eikonal}) as%
\[
U^{2}+\left(  \frac{p-y}{z}+q-p\right)  U+\left(  U_{0}\right)  ^{2}=0
\]
and solving for $U$ we get%
\begin{equation}
U^{\pm}(y,z)=-\frac{1}{2}\left(  \frac{p-y}{z}+q-p\right)  \pm\frac{1}{2}%
\sqrt{\left(  \frac{p-y}{z}+q-p\right)  ^{2}-4\left(  U_{0}\right)  ^{2}}
\label{U+-}%
\end{equation}
where%
\begin{equation}
U_{0}(z)=\sqrt{\frac{pq(1-z)}{z}}. \label{Uo}%
\end{equation}

\begin{figure}[ptb]
\begin{center}
\rotatebox{270} {\resizebox{12cm}{!}{\includegraphics{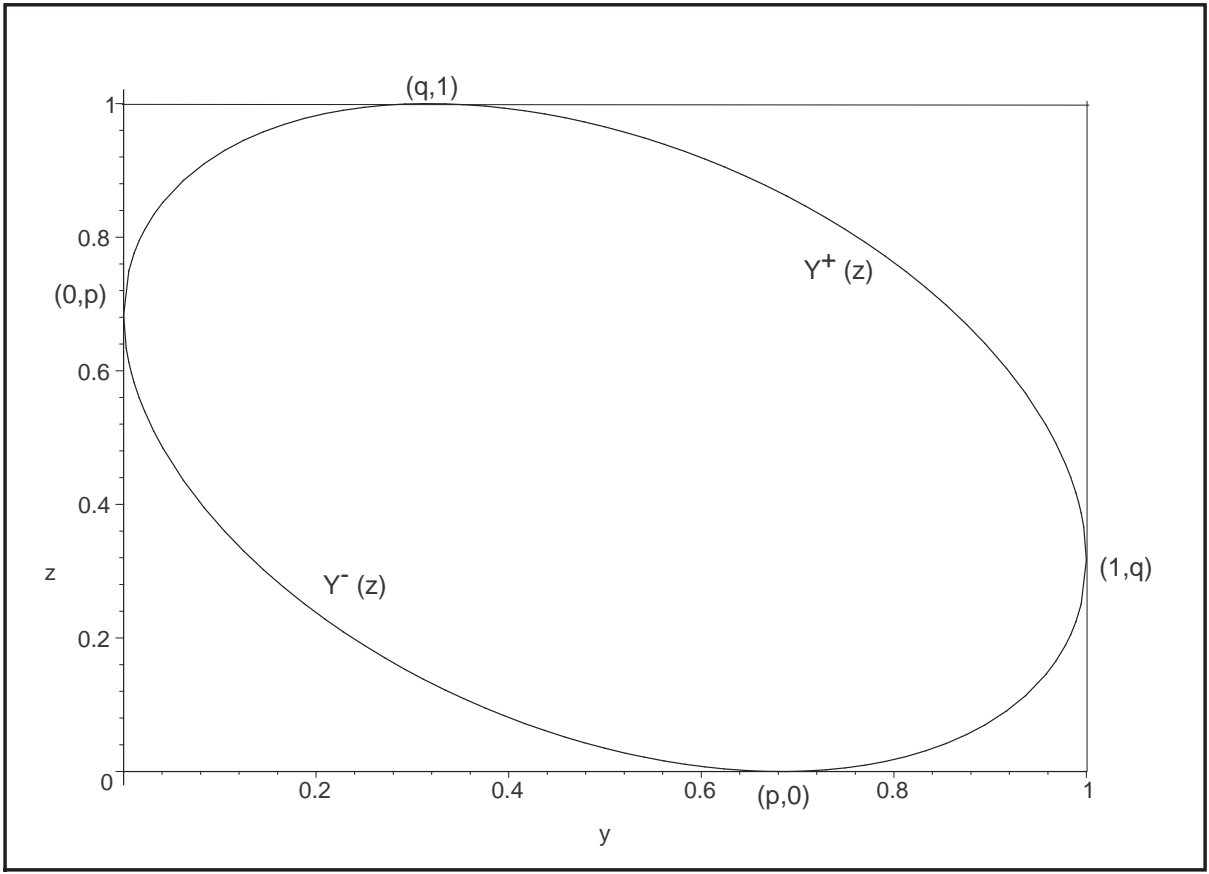}}}
\end{center}
\caption{A sketch of the ellipse \textbf{$\mathbb{E}$} and the curves $Y^{\pm
}(z)$.}%
\label{elipse}%
\end{figure}

The discriminant in (\ref{U+-}) vanishes if $\frac{p-y}{z}+q-p=\pm2U_{0},$
which is equivalent to $y=Y^{\pm},$ with%
\begin{equation}
Y^{\pm}(z)=p+\left(  q-p\right)  z\pm2zU_{0}. \label{Y+-}%
\end{equation}
Rewriting the equation $\left(  \frac{p-y}{z}+q-p\right)  ^{2}-4\left(
U_{0}\right)  ^{2}=0$ as
\begin{equation}
\left(  y-\frac{1}{2}\right)  ^{2}+\left(  z-\frac{1}{2}\right)
^{2}+2(p-q)\left(  y-\frac{1}{2}\right)  \left(  z-\frac{1}{2}\right)  -pq=0
\label{eqE}%
\end{equation}
we can see that
\begin{equation}
\left\{  (y,z)\in\lbrack0,1]\times\lbrack0,1]:y=Y^{\pm}(z)\right\}
=\mathbf{\mathbb{E}} \label{E}%
\end{equation}
where \textbf{$\mathbb{E}$} is an ellipse centered at $\left(  \frac{1}%
{2},\frac{1}{2}\right)  $ (see Figure \ref{elipse}). After rotation by
$\pm\frac{\pi}{4}$ and translation to the origin, \textbf{$\mathbb{E}$}
reduces to one of the canonical forms%
\begin{align*}
\frac{\overline{y}^{2}}{q/2}+\frac{\overline{z}^{2}}{p/2}  &  =1,\qquad p<q\\
\overline{y}^{2}+\overline{z}^{2}  &  =\frac{1}{4},\qquad p=\frac{1}{2}=q\\
\frac{\overline{y}^{2}}{p/2}+\frac{\overline{z}^{2}}{q/2}  &  =1,\qquad p>q.
\end{align*}

The ellipse \textbf{$\mathbb{E}$} is contained in the square $[0,1]\times
\lbrack0,1]$ and intersects the $y$ and $z$ axis at the points
$(0,p),(1,q),(p,0)$ and $(q,1).$ Its left side between the points $(q,1)$ and
$(p,0)$ coincides with the curve $Y^{-}(z)$ and its right side with the curve
$Y^{+}(z).$

For points $(y,z)$ located outside \textbf{$\mathbb{E}$}$\mathbf{,}$ $U^{\pm
}(y,z)$ are real and for points $(y,z)$ inside \textbf{$\mathbb{E}$} $U^{\pm
}(y,z)$ are complex conjugates. When $y=Y^{\pm}$ the two values $U^{\pm}$
coalesce and we have%
\begin{equation}
U^{+}(Y^{+},z)=U_{0}=U^{-}(Y^{+},z),\quad U^{+}(Y^{-},z)=-U_{0}=U^{-}%
(Y^{-},z). \label{U(Y+)}%
\end{equation}

Writing the function $L(y,z)$ in terms of $U_{0}$ we have from (\ref{L1})%
\begin{equation}
L(y,z)=B(y)\sqrt{\frac{(U-p)(U+q)}{z\left(  U^{2}-U_{0}^{2}\right)  }}.
\label{L}%
\end{equation}
In the rest of this paper we shall use the following notation%
\begin{equation}
\psi^{\pm}(y,z)=\ln\left[  \left(  U^{\pm}\right)  ^{z-1}\left(  U^{\pm
}-p\right)  ^{1-y}\left(  U^{\pm}+q\right)  ^{y}\right]  \label{psi+-}%
\end{equation}
and%
\begin{equation}
L^{\pm}(y,z)=\sqrt{\frac{(U^{\pm}-p)(U^{\pm}+q)}{z\left[  \left(  U^{\pm
}\right)  ^{2}-U_{0}^{2}\right]  }}. \label{L+-}%
\end{equation}
Hence, we write
\begin{align}
K_{n}(x)  &  \sim\varepsilon^{\nu}B^{-}(y)\exp\left[  \varepsilon^{-1}\psi
^{-}(y,z)+\varepsilon^{-1}A^{-}(y)\right]  L^{-}(y,z)\label{F+-}\\
&  +\varepsilon^{\nu}B^{+}(y)\exp\left[  \varepsilon^{-1}\psi^{+}%
(y,z)+\varepsilon^{-1}A^{+}(y)\right]  L^{+}(y,z)\nonumber
\end{align}
where $B^{\pm}(y)$ and $A^{\pm}(y)$ are functions to be determined. From
(\ref{L+-}) we see that $L^{\pm}(y,z)$ are singular when $U^{\pm}=U_{0},$
i.e., for $y=Y^{\pm},$ and also for $z=0.$ Therefore, we need to find
asymptotic solutions valid in those regions.

\section{The boundary $n=0$ (Region I)}

For $n=O(1)$ we have from (\ref{Kndef}) that $K_{n}(x)=O\left(  N^{n}\right)
$ as $N\rightarrow\infty.$ Thus, we introduce the function $R_{n}^{(1)}\left(
y\right)  $ and consider solutions of (\ref{recurrence}) which have the
asymptotic form%
\begin{equation}
K_{n}(x)=N^{n}R_{n}^{(1)}\left(  \frac{x}{N}\right)  . \label{Knz=0}%
\end{equation}
Using (\ref{Knz=0}) in (\ref{recurrence}) and expanding in powers of $N$
gives, to leading order%
\[
(n+1)R_{n+1}^{(1)}+(p-y)R_{n}^{(1)}=0.
\]
Solving the recursion above with the initial condition (\ref{condn=0}) we get%
\[
R_{n}^{(1)}(y)=\frac{\left(  y-p\right)  ^{n}}{n!}%
\]
and hence%
\begin{equation}
K_{n}(x)\sim K_{n}^{(1)}(y)=N^{n}\frac{\left(  y-p\right)  ^{n}}{n!},\qquad
n=O(1). \label{K1}%
\end{equation}

\subsection{The corner layer at $(p,0)$ (Region II)}

We shall now find an asymptotic solution in the neighborhood of the point
$(p,0).$ We introduce the stretched variable $\eta$ and the function
$R_{n}^{(2)}(\eta)$ defined by%
\begin{align}
y &  =p+\eta\sqrt{2pq\varepsilon}\quad\eta=O(1),\label{eta}\\
\quad\quad K_{n}(x) &  =\frac{\varepsilon^{-n/2}}{n!}\left(  \frac{pq}%
{2}\right)  ^{n/2}R_{n}^{(2)}\left(  \frac{\varepsilon x-p}{\sqrt
{2pq\varepsilon}}\right)  .\nonumber
\end{align}
Using (\ref{eta}) in (\ref{recurrence}) yields, to leading order, the equation%
\[
R_{n+1}^{(2)}-2\eta R_{n}^{(2)}+2nR_{n-1}^{(2)}=0
\]
which we recognize as the recurrence relation for the Hermite polynomials.
Thus,%
\begin{equation}
K_{n}(x)\sim K_{n}^{\left(  2\right)  }(\eta)=\frac{\varepsilon^{-n/2}}%
{n!}\left(  \frac{pq}{2}\right)  ^{n/2}H_{n}\left(  \eta\right)  \label{Hn}%
\end{equation}
for $n=O(1)$ and $y-p=O\left(  \varepsilon^{1/2}\right)  ,$ where
$H_{n}\left(  \eta\right)  $ is the Hermite polynomial of degree $n$.

\section{The lower corners (Regions III and IV)}

Setting $n=z/\varepsilon$ in (\ref{K1}) and letting $\varepsilon\rightarrow0$
we obtain
\begin{equation}
K_{n}^{(1)}(y)\sim\varepsilon^{1/2}\frac{1}{\sqrt{2\pi z}}\exp\left\{
\varepsilon^{-1}\left[  1-\ln(z)+\ln(y-p)\right]  z\right\}  \label{K1large}%
\end{equation}
where we have used Stirling's formula \cite{MR94b:00012}%
\begin{equation}
\Gamma(x)\sim\sqrt{\frac{2\pi}{x}}x^{x}e^{-x},\quad x\rightarrow\infty.
\label{stirling}%
\end{equation}
From (\ref{psi+-})-(\ref{L+-}) we have as $z\rightarrow0$%
\begin{equation}
\psi^{-}(y,z)\sim\left\{
\begin{array}
[c]{c}%
\left[  1-\ln(z)+\ln(y-p)\right]  z,\quad y<p\\
(1-y)\ln\left(  \frac{1-y}{q}\right)  +y\ln\left(  \frac{y}{p}\right)
+\ln\left(  \frac{pq}{y-p}\right)  z,\quad y>p
\end{array}
\right.  \label{psimz=0}%
\end{equation}%
\begin{equation}
L^{-}(y,z)=\left\{
\begin{array}
[c]{c}%
z^{-1/2}+O(z),\quad y<p\\
\frac{\sqrt{(y-1)y}}{y-p}+O(z),\quad y>p
\end{array}
\right.  \label{L-z=0}%
\end{equation}%
\begin{equation}
\psi^{+}(y,z)\sim\left\{
\begin{array}
[c]{c}%
(1-y)\ln\left(  \frac{1-y}{q}\right)  +y\ln\left(  \frac{y}{p}\right)
+\ln\left(  \frac{pq}{y-p}\right)  z,\quad y<p\\
\left[  1-\ln(z)+\ln(y-p)\right]  z,\quad y>p
\end{array}
\right.  \label{psi+z=0}%
\end{equation}
and%
\begin{equation}
L^{+}(y,z)=\left\{
\begin{array}
[c]{c}%
\frac{\sqrt{(y-1)y}}{p-y}+O(z),\quad y<p\\
z^{-1/2}+O(z),\quad y>p
\end{array}
\right.  . \label{L+z=0}%
\end{equation}
Matching (\ref{psimz=0})-(\ref{L+z=0}) and (\ref{K1large}) we conclude that
\begin{equation}
K_{n}(x)\sim\left\{
\begin{array}
[c]{c}%
K^{-}(y,z),\quad0<y<Y^{-}(z)\\
K^{+}(y,z),\quad Y^{+}(z)<y<1
\end{array}
\right.  \label{K3,4}%
\end{equation}
with%
\begin{equation}
K^{-}(y,z)=\varepsilon^{1/2}\frac{1}{\sqrt{2\pi}}\exp\left[  \varepsilon
^{-1}\psi^{-}(y,z)\right]  L^{-}(y,z) \label{K-}%
\end{equation}%
\begin{equation}
K^{+}(y,z)=\varepsilon^{1/2}\frac{1}{\sqrt{2\pi}}\exp\left[  \varepsilon
^{-1}\psi^{+}(y,z)\right]  L^{+}(y,z). \label{K+}%
\end{equation}

\begin{remark}
In the remainder of the paper, we will find asymptotic formulas only in the
region $0\leq y\leq Y^{+}(z),$ $0\leq z\leq1.$ The corresponding results for
$\ Y^{+}(z)\leq y\leq1$ can be obtained by using the symmetry formula
(\ref{symetry}) and noting that under the transformations $y\rightarrow1-y,$
$p\leftrightarrow q,$ $K_{n}\rightarrow\left(  -1\right)  ^{n}K_{n},$ we
obtain
\begin{align*}
U^{-}  &  \rightarrow U^{+}\\
\psi^{-}  &  \rightarrow\psi^{+}+z\pi\mathrm{i,\quad}L^{-}\rightarrow L^{+}\\
K^{-}  &  \rightarrow K^{+}.
\end{align*}

\end{remark}

\section{The boundary $x=0$}

We shall now consider the case $x=O(1).$

\begin{lemma}
Let $x=m,$ with $m$ an integer, $m\ll N.$

\begin{enumerate}
\item
\begin{equation}
K_{n}(m)=\left(  -p\right)  ^{n}\binom{N}{n}\left(  1-\frac{n}{p}%
N^{-1}\right)  ^{m},\quad m=0,1. \label{x=O(1)1}%
\end{equation}

\item If $n=O(N),$ then%
\begin{equation}
K_{n}(m)\sim\left(  -p\right)  ^{n}\binom{N}{n}\left(  1-\frac{n}{p}%
N^{-1}\right)  ^{m},\quad N\rightarrow\infty,\quad m\geq2. \label{x=O(1)2}%
\end{equation}

\end{enumerate}
\end{lemma}

\begin{proof}
From (\ref{Kndef}) we have for $x=m$ integer%
\begin{align}
K_{n}(m)  &  =\sum_{k=0}^{m}\binom{m}{k}\binom{N-m}{n-k}q^{k}\left(
-p\right)  ^{n-k}\label{x=O(1)sum}\\
&  =\left(  -p\right)  ^{n}\binom{N}{n}\sum_{k=0}^{m}\binom{m}{k}\frac
{\binom{N-m}{n-k}}{\binom{N}{n}}\left(  -\frac{q}{p}\right)  ^{k}\nonumber
\end{align}
and (\ref{x=O(1)1}) follows for $m=0,1.$

Setting $n=zN$ and using (\ref{stirling}) we get%
\[
\frac{\binom{N-m}{Nz-k}}{\binom{N}{Nz}}\sim\left(  1-z\right)  ^{m}\left(
\frac{z}{1-z}\right)  ^{k},\quad N\rightarrow\infty.
\]
Using the above in (\ref{x=O(1)sum}) we have%
\begin{align*}
K_{n}(m)  &  \sim\left(  -p\right)  ^{n}\binom{N}{n}\left(  1-z\right)
^{m}\sum_{k=0}^{m}\binom{m}{k}\left(  -\frac{z}{1-z}\frac{q}{p}\right)  ^{k}\\
&  =\left(  -p\right)  ^{n}\binom{N}{n}\left(  1-z\right)  ^{m}\left[
\frac{p-z}{p\left(  1-z\right)  }\right]  ^{m}%
\end{align*}
and (\ref{x=O(1)2}) follows.
\end{proof}

Using (\ref{stirling}) we have, as $N\rightarrow\infty$%
\begin{equation}
\left(  -p\right)  ^{n}\binom{N}{n}\sim\frac{\varepsilon^{1/2}}{\sqrt{2\pi
}\sqrt{z\left(  1-z\right)  }}\exp\left[  \varepsilon^{-1}\phi_{0}(z)\right]
\label{binomial}%
\end{equation}
where%
\begin{equation}
\phi_{0}(z)=(z-1)\ln(1-z)-z\ln(z)+z\ln(-p). \label{phi0}%
\end{equation}
From (\ref{K-}) we get, as $y\rightarrow0$%
\begin{equation}
K^{-}(y,z)\sim\frac{\varepsilon^{1/2}}{\sqrt{2\pi}\sqrt{z\left(  1-z\right)
}}\exp\left[  \frac{\phi_{0}(z)+\ln\left(  \frac{p-z}{p}\right)
y}{\varepsilon}\right]  ,\quad0<z<p. \label{K-y=0}%
\end{equation}
Hence, as $y\rightarrow0,$ $K^{-}(y,z)$ satisfies the boundary condition
(\ref{condx=0}) for $0<z<p$.

\subsection{The boundary layer at $x=0,$ $p<z<1$ (Region V)}

Taking (\ref{x=O(1)2}) into account, we define the function $R_{n}^{(5)}(x)$
by%
\begin{equation}
K_{n}(x)=\binom{N}{n}\left(  -p\right)  ^{n}R_{n}^{(5)}(x). \label{R12}%
\end{equation}
Using (\ref{R12}) in (\ref{recurrence}) yields
\begin{equation}
p(n-N)R_{n+1}^{(5)}+\left[  pN+\left(  q-p\right)  n-x\right]  R_{n}%
^{(5)}-nqR_{n-1}^{(5)}=0. \label{recurr12}%
\end{equation}
Writing (\ref{recurr12}) in terms of $z,\varepsilon$ and the function
$G^{(5)}(x,z)$ defined by%
\begin{equation}
R_{n}^{(5)}(x)=G^{(5)}(x,\varepsilon n) \label{G}%
\end{equation}
we have%
\begin{equation}
p(1-z)G^{(5)}(x,z+\varepsilon)+\left[  z\left(  p-q\right)  +x\varepsilon
-p\right]  G^{(5)}(x,z)+zqG^{(5)}(x,z-\varepsilon)=0. \label{de2}%
\end{equation}
Using the WKB anszat%
\begin{equation}
G^{(5)}(x,z)\sim\varepsilon^{\tau}\exp\left[  \frac{1}{\varepsilon}%
\phi(x,z)\right]  W(x,z) \label{WKBG}%
\end{equation}
in (\ref{de2}), we obtain the equations%
\begin{equation}
\left[  \exp\left(  \phi_{z}\right)  -1\right]  \left[  p(z-1)\exp\left(
\phi_{z}\right)  +zq\right]  W=0 \label{phi}%
\end{equation}
and%
\begin{gather}
\left[  zq+p(z-1)\exp\left(  2\phi_{z}\right)  \right]  W_{z}\label{W}\\
+\left\{  \frac{1}{2}\left[  p(z-1)\exp\left(  2\phi_{z}\right)  -zq\right]
\phi_{zz}-x\exp\left(  \phi_{z}\right)  \right\}  W=0.\nonumber
\end{gather}
Solving (\ref{phi})-(\ref{W}) we get%
\begin{gather}
G^{(5)}(x,z)\sim\varepsilon^{\tau_{1}}A^{(5)}(x)\exp\left[  \varepsilon
^{-1}\varphi_{1}(x)\right]  \left(  z-p\right)  ^{x}\label{G1}\\
+\varepsilon^{\tau_{2}}B^{(5)}(x)\exp\left[  \varepsilon^{-1}\phi
_{2}(x,z)\right]  \left(  z-p\right)  ^{-x-1}\sqrt{z\left(  1-z\right)
}\nonumber
\end{gather}
where%
\begin{equation}
\phi_{2}(x,z)=\left(  1-z\right)  \ln\left(  1-z\right)  +z\ln(z)+(z-1)\ln
(q)-z\ln(p)+\varphi_{2}(x) \label{phi2}%
\end{equation}
and the coefficients $\tau_{1},\tau_{2}$ and the functions $A^{(5)}%
(x),B^{(5)}(x),\varphi_{1}(x),\varphi_{2}(x)$ are to be determined.

Using (\ref{binomial}) we obtain%
\begin{gather}
K_{n}(x)\sim K^{(5)}\left(  x,z\right)  =\frac{\varepsilon^{1/2+\tau_{1}}%
}{\sqrt{2\pi}\sqrt{z\left(  1-z\right)  }}\exp\left[  \frac{\phi
_{0}(z)+\varphi_{1}(x)}{\varepsilon}\right]  A(x)\left(  z-p\right)
^{x}\label{K123}\\
+\frac{\varepsilon^{1/2+\tau_{2}}}{\sqrt{2\pi}}\exp\left[  \frac
{(z-1)\ln(q)+\pi\mathrm{i}z+\varphi_{2}(x)}{\varepsilon}\right]  B(x)\left(
z-p\right)  ^{-x-1}.\nonumber
\end{gather}

\subsection{The corner layer at $(0,p)$ (Region VI)}

For $x\approx0$ and $n\approx Np,$ we scale $n$ as
\begin{equation}
n=Np-u\sqrt{pqN},\quad u=O(1) \label{u}%
\end{equation}
and introduce the function $G^{(6)}(x,u)$ defined by%
\begin{equation}
K_{n}(x)=\binom{N}{n}\left(  -p\right)  ^{n}G^{(6)}\left(  x,\frac{Np-n}%
{\sqrt{pqN}}\right)  . \label{R13}%
\end{equation}
Using (\ref{u})-(\ref{R13}) in (\ref{recurr12}) we get, to leading order%
\begin{equation}
\frac{\partial^{2}G^{(6)}}{\partial u^{2}}-\frac{\partial G^{(6)}}{\partial
u}+xG^{(6)}=0. \label{R13eq}%
\end{equation}
Solving (\ref{R13eq}) we obtain%
\begin{equation}
G^{(6)}(x,u)=\exp\left(  \frac{u^{2}}{4}\right)  \left[  A^{(6)}%
(x)\mathrm{D}_{x}(u)+B^{(6)}(x)\mathrm{D}_{x}(-u)\right]  \label{D}%
\end{equation}
where \textrm{$D$}$_{x}(u)$ is the parabolic cylinder function and
$A^{(6)}(x),$ $B^{(6)}(x)$ are functions to be determined. Since%
\begin{equation}
u=\frac{Np-n}{\sqrt{pqN}}=\frac{p-z}{\sqrt{pq\varepsilon}}, \label{u2}%
\end{equation}
we note that the limit $u\rightarrow\infty$ corresponds to the matching
between regions VI and III, while the limit $u\rightarrow-\infty$ corresponds
to the matching between regions VI and V. As $z\rightarrow p$ we have%
\begin{gather}
\left(  -p\right)  ^{n}\binom{N}{n}\sim\frac{\sqrt{\varepsilon}}{\sqrt{2\pi
pq}}\exp\left(  -\frac{u^{2}}{2}\right) \label{binomial1}\\
\times\exp\left\{  \frac{\pi\mathrm{i}p-q\ln\left(  q\right)  }{\varepsilon
}-\frac{u\sqrt{pq}\left[  \pi\mathrm{i}+\ln\left(  q\right)  \right]  }%
{\sqrt{\varepsilon}}\right\}  .\nonumber
\end{gather}
Thus,%
\begin{gather}
K_{n}(x)\sim K^{(6)}(x,u)=\frac{\sqrt{\varepsilon}}{\sqrt{2\pi pq}}\left[
A^{(6)}(x)\mathrm{D}_{x}(u)+B^{(6)}(x)\mathrm{D}_{x}(-u)\right] \label{K132}\\
\times\exp\left\{  \frac{\pi\mathrm{i}p-q\ln\left(  q\right)  }{\varepsilon
}-\frac{u\sqrt{pq}\left[  \pi\mathrm{i}+\ln\left(  q\right)  \right]  }%
{\sqrt{\varepsilon}}-\frac{u^{2}}{4}\right\}  .\nonumber
\end{gather}

Using (\ref{u2}) in (\ref{K-y=0}) yields%
\begin{gather}
K^{-}\left(  y,z\right)  \sim\frac{\sqrt{\varepsilon}}{\sqrt{2\pi pq}}\left[
\sqrt{\frac{q\varepsilon}{p}}u\right]  ^{x}\label{u1}\\
\times\exp\left[  \frac{\pi\mathrm{i}p-q\ln\left(  q\right)  }{\varepsilon
}-\frac{u\sqrt{pq}\left[  \pi\mathrm{i}+\ln\left(  q\right)  \right]  }%
{\sqrt{\varepsilon}}-\frac{u^{2}}{2}\right]  .\nonumber
\end{gather}

Using the well known asymptotic approximation \cite{MR97k:01072}%
\begin{equation}
\mathrm{D}_{x}(u)\sim\exp\left(  -\frac{u^{2}}{4}\right)  u^{x},\quad
u\rightarrow\infty\label{Dularge}%
\end{equation}
in (\ref{K132}) and comparing with (\ref{u1}) we conclude that
\[
A^{(6)}(x)=\left[  \sqrt{\frac{q\varepsilon}{p}}\right]  ^{x}%
\]
and $B^{(6)}(x)=0.$ Therefore,%
\begin{gather}
K^{(6)}(x,u)=\frac{\varepsilon^{1/2}}{\sqrt{2\pi pq}}\left[  \sqrt
{\frac{q\varepsilon}{p}}\right]  ^{x}\mathrm{D}_{x}(u)\label{K13}\\
\times\exp\left[  \frac{\pi\mathrm{i}p-q\ln\left(  q\right)  }{\varepsilon
}-\frac{u\sqrt{pq}\left[  \pi\mathrm{i}+\ln\left(  q\right)  \right]  }%
{\sqrt{\varepsilon}}-\frac{u^{2}}{4}\right] \nonumber
\end{gather}
for $x=O(1),$ $z-p=O\left(  \varepsilon^{1/2}\right)  .$

Using the formula \cite{MR97k:01072}%
\begin{equation}
\mathrm{D}_{x}(-u)\sim\exp\left(  -\frac{u^{2}}{4}\right)  u^{x}\cos(\pi
x)-\sqrt{\frac{2}{\pi}}x\Gamma(x)\sin(\pi x)u^{-x-1}\exp\left(  \frac{u^{2}%
}{4}\right)  ,\quad u\rightarrow\infty\label{Duneg}%
\end{equation}
in (\ref{K13}) we have%
\begin{gather}
K^{(6)}(x,-u)\sim\frac{\varepsilon^{1/2}}{\sqrt{2\pi pq}}\left[  u\sqrt
{\frac{q\varepsilon}{p}}\right]  ^{x}\cos(\pi x)\nonumber\\
\times\exp\left[  \frac{\pi\mathrm{i}p-q\ln\left(  q\right)  }{\varepsilon
}-\frac{u\sqrt{pq}\left[  \pi\mathrm{i}+\ln\left(  q\right)  \right]  }%
{\sqrt{\varepsilon}}-\frac{u^{2}}{2}\right] \label{K135}\\
-\frac{1}{u}\frac{\varepsilon^{1/2}}{\pi\sqrt{pq}}\left[  \frac{1}{u}%
\sqrt{\frac{q\varepsilon}{p}}\right]  ^{x}x\Gamma(x)\sin(\pi x)\nonumber\\
\times\exp\left[  \frac{\pi\mathrm{i}p-q\ln\left(  q\right)  }{\varepsilon
}-\frac{u\sqrt{pq}\left[  \pi\mathrm{i}+\ln\left(  q\right)  \right]  }%
{\sqrt{\varepsilon}}\right]  ,\quad u\rightarrow\infty.\nonumber
\end{gather}

Using (\ref{u2}) in (\ref{K123}) gives, as $z\downarrow p$%
\begin{gather}
K^{(5)}\left(  x,z\right)  \sim\frac{\varepsilon^{1/2+\tau_{1}}}{\sqrt{2\pi
pq}}A^{(5)}(x)\left(  u\sqrt{pq\varepsilon}\right)  ^{x}\nonumber\\
\times\exp\left[  \frac{\pi\mathrm{i}p-q\ln\left(  q\right)  +\varphi_{1}%
(x)}{\varepsilon}-\frac{u\sqrt{pq}\left[  \pi\mathrm{i}+\ln\left(  q\right)
\right]  }{\sqrt{\varepsilon}}-\frac{u^{2}}{2}\right] \label{K125}\\
+\exp\left[  \frac{\pi\mathrm{i}p-q\ln\left(  q\right)  +\varphi_{2}%
(x)}{\varepsilon}-\frac{u\sqrt{pq}\left[  \pi\mathrm{i}+\ln\left(  q\right)
\right]  }{\sqrt{\varepsilon}}\right] \nonumber\\
\times\frac{1}{u}\frac{\varepsilon^{\tau_{2}}}{\sqrt{2\pi pq}}B^{(5)}%
(x)\left(  u\sqrt{pq\varepsilon}\right)  ^{-x}.\nonumber
\end{gather}
Matching (\ref{K135}) and (\ref{K125}) yields $\tau_{1}=0,$ $\varphi
_{1}(x)=0,$ $\tau_{2}=\frac{1}{2},$ $\varphi_{2}(x)=0$ and
\begin{align*}
A^{(5)}(x)  &  =p^{-x}\cos(\pi x)\\
B^{(5)}(x)  &  =-\sqrt{\frac{2}{\pi}}x\Gamma(x)\sin(\pi x)\left(
q\varepsilon\right)  ^{x}.
\end{align*}
Thus,%
\begin{gather}
K^{(5)}\left(  x,z\right)  =\frac{\varepsilon^{1/2}}{\sqrt{2\pi}\sqrt{z\left(
1-z\right)  }}\cos(\pi x)\left(  \frac{z-p}{p}\right)  ^{x}\exp\left[
\frac{\phi_{0}(z)}{\varepsilon}\right] \label{K12}\\
-\frac{\varepsilon}{\pi}\frac{x}{z-p}\Gamma(x)\sin(\pi x)\left(
\frac{q\varepsilon}{z-p}\right)  ^{x}\exp\left[  \frac{(z-1)\ln(q)+\pi
\mathrm{i}z}{\varepsilon}\right] \nonumber
\end{gather}
for $x=O(1)$ and $p<z<1.$

From (\ref{K12}) we have, as $x\rightarrow\infty$%
\begin{gather}
K^{(5)}\left(  x,z\right)  \sim\frac{\varepsilon^{1/2}}{\sqrt{2\pi}%
\sqrt{z\left(  1-z\right)  }}\cos(\pi x)\left(  \frac{z-p}{p}\right)  ^{x}%
\exp\left[  \frac{\phi_{0}(z)}{\varepsilon}\right] \label{K12xlarge}\\
-\varepsilon\sqrt{\frac{2}{\pi}}\frac{\sqrt{x}}{z-p}\sin(\pi x)\left(
\frac{q\varepsilon x}{z-p}\right)  ^{x}\exp\left[  \frac{(z-1)\ln
(q)+\pi\mathrm{i}z}{\varepsilon}-x\right]  .\nonumber
\end{gather}
In terms of $y=x\varepsilon,$ (\ref{K12xlarge}) reads%
\begin{gather}
K^{(5)}\left(  x,z\right)  \sim\frac{\varepsilon^{1/2}}{\sqrt{2\pi}%
\sqrt{z\left(  1-z\right)  }}\cos\left(  \frac{\pi y}{\varepsilon}\right)
\exp\left[  \frac{\phi_{0}(z)+y\ln\left(  \frac{z-p}{p}\right)  }{\varepsilon
}\right] \label{K12y}\\
-\varepsilon^{1/2}\sqrt{\frac{2}{\pi}}\sin\left(  \frac{\pi y}{\varepsilon
}\right)  \frac{\sqrt{y}}{z-p}\exp\left[  \frac{(z-1)\ln(q)+\pi\mathrm{i}%
z+y\ln\left(  \frac{qy}{z-p}\right)  -y}{\varepsilon}\right]  .\nonumber
\end{gather}

\section{The left upper corner (Region VII)}

We now consider the region $0\ll y<Y^{-}(z),$ $p<z<1.$ From (\ref{psi+-}%
)-(\ref{L+-}) we have for $y\rightarrow0,\quad z>p$%
\begin{equation}
\psi^{+}(y,z)\sim\phi_{0}(z)+y\ln\left(  \frac{z-p}{p}\right)  -\pi
\mathrm{i}y, \label{psipy=0}%
\end{equation}%
\begin{equation}
\psi^{-}(y,z)\sim(z-1)\ln(q)+\pi\mathrm{i}z+y\ln\left(  \frac{qy}{z-p}\right)
-y-\pi\mathrm{i}y, \label{psimy=0}%
\end{equation}
and%
\begin{equation}
L^{+}(y,z)\sim\frac{1}{\sqrt{z\left(  1-z\right)  }},\quad L^{-}(y,z)\sim
\frac{\sqrt{y}}{z-p}\mathrm{i.} \label{Ly=0}%
\end{equation}
Using (\ref{psipy=0})-(\ref{Ly=0}) in (\ref{F+-}) and matching with
(\ref{K12y}), we conclude that%
\begin{equation}
K_{n}(x)\sim K^{\left(  7\right)  }(y,z)=\exp\left(  \frac{\pi\mathrm{i}%
y}{\varepsilon}\right)  \left[  \cos\left(  \frac{\pi y}{\varepsilon}\right)
K^{+}(y,z)+2\mathrm{i}\sin\left(  \frac{\pi y}{\varepsilon}\right)
K^{-}(y,z)\right]  , \label{leftcorner}%
\end{equation}
which can be written as%
\begin{equation}
K^{\left(  7\right)  }(y,z)=\frac{1}{2}\left[  \exp\left(  \frac
{2\pi\mathrm{i}y}{\varepsilon}\right)  +1\right]  K^{+}(y,z)+\left[
\exp\left(  \frac{2\pi\mathrm{i}y}{\varepsilon}\right)  -1\right]  K^{-}(y,z).
\label{left1}%
\end{equation}

\section{The transition layer $y=Y^{-}(z)$}

As we noted before, the functions $L^{\pm}(y,z)$ are infinite on the curves
$y=Y^{\pm}(z).$ Hence, we need to find transition layer solutions there.

\subsection{The lower part $0<z<p$ (Region VIII)}

We introduce the stretched variable $\beta,$ defined by%
\begin{equation}
y=Y^{-}(z)-\beta\varepsilon^{2/3}\quad\beta=O(1),\quad0<z<p. \label{beta}%
\end{equation}
Using (\ref{beta}) in (\ref{psi+-})-(\ref{L+-}) and expanding in powers of
$\varepsilon,$ with $\beta>0,$ we get%
\begin{equation}
\psi^{-}(y,z)\sim\psi_{0}(z)+\ln\left(  \frac{U_{0}+p}{U_{0}-q}\right)
\beta\varepsilon^{2/3}-\frac{2}{3}\sqrt{\frac{U_{0}}{z}}\frac{1}{\left(
U_{0}+p\right)  \left(  U_{0}-q\right)  }\beta^{3/2}\varepsilon\label{psi-Y-}%
\end{equation}
where%
\begin{equation}
\psi_{0}(z)=z\pi\mathrm{i}+(z-1)\ln\left(  U_{0}\right)  +Y^{-}(z)\ln\left(
U_{0}-q\right)  +\left[  1-Y^{-}(z)\right]  \ln\left(  U_{0}+p\right)
\label{psi0-}%
\end{equation}
and%
\begin{equation}
L^{-}(y,z)\sim\frac{1}{\sqrt{2}}z^{-1/4}\sqrt{\left(  U_{0}+p\right)  \left(
U_{0}-q\right)  }\left(  U_{0}\right)  ^{-3/4}\beta^{-1/4}\varepsilon^{-1/6}.
\label{L-Y-}%
\end{equation}
Hence,%
\begin{gather}
K^{-}(y,z)\sim\varepsilon^{1/3}\frac{1}{2\sqrt{\pi}}z^{-1/4}\sqrt{\left(
U_{0}+p\right)  \left(  U_{0}-q\right)  }\left(  U_{0}\right)  ^{-3/4}%
\beta^{-1/4}\label{K-1}\\
\times\exp\left[  \varepsilon^{-1}\psi_{0}(z)+\ln\left(  \frac{U_{0}+p}%
{U_{0}-q}\right)  \beta\varepsilon^{-1/3}-\frac{2}{3}\sqrt{\frac{U_{0}}{z}%
}\frac{1}{\left(  U_{0}+p\right)  \left(  U_{0}-q\right)  }\beta^{3/2}\right]
\nonumber
\end{gather}
for $y\uparrow Y^{-}(z).$ Thus, we introduce the function $G^{(8)}(\beta,z)$
and consider solutions of the form%
\begin{equation}
K_{n}(x)=\varepsilon^{\nu_{8}}\exp\left[  \psi_{0}(z)\varepsilon^{-1}%
+\ln\left(  \frac{U_{0}+p}{U_{0}-q}\right)  \left(  Y^{-}-\varepsilon
x\right)  \varepsilon^{-1}\right]  G^{(8)}\left(  \frac{Y^{-}-\varepsilon
x}{\varepsilon^{2/3}},\varepsilon n\right)  . \label{G8}%
\end{equation}
Using (\ref{G8}) in (\ref{recurrence}) and expanding in powers of
$\varepsilon$ we obtain%
\begin{equation}
G_{\beta\beta}^{(8)}=-\frac{2\beta}{\left[  \left(  Y^{-}\right)  ^{\prime
}\right]  ^{2}\left[  \left(  p-q\right)  z+Y^{-}-p\right]  }G^{(8)}%
=\beta\Theta^{2}G^{(8)} \label{Gbeta}%
\end{equation}
where
\begin{equation}
\Theta(z)=\sqrt{\frac{U_{0}}{z}}\frac{1}{\left(  U_{0}+p\right)  \left(
U_{0}-q\right)  }. \label{Theta-}%
\end{equation}
Solving (\ref{Gbeta}) we get%
\begin{equation}
G^{(8)}(\beta,z)=A^{(8)}(z)\mathrm{Ai}\left[  \Theta^{2/3}\beta\right]
+B^{(8)}(z)\mathrm{Bi}\left[  \Theta^{2/3}\beta\right]  \label{G3}%
\end{equation}
where $A^{(8)}(z),B^{(8)}(z)$ are functions to be determined and
$\mathrm{Ai}\left(  \cdot\right)  ,\mathrm{Bi}\left(  \cdot\right)  $ are the
Airy functions. Using the formulas \cite{MR94b:00012}%
\begin{align}
\mathrm{Ai}\left(  x\right)   &  \sim\frac{1}{2\sqrt{\pi}}x^{-1/4}\exp\left(
-\frac{2}{3}x^{3/2}\right)  ,\quad x\rightarrow\infty\label{Ailarge}\\
\mathrm{Bi}\left(  x\right)   &  \sim\frac{1}{\sqrt{\pi}}x^{-1/4}\exp\left(
\frac{2}{3}x^{3/2}\right)  ,\quad x\rightarrow\infty\nonumber
\end{align}
in (\ref{G3}) and matching with (\ref{K-1}) we conclude that $\nu_{8}=1/3,$
$B^{(8)}(z)=0$ and%
\begin{equation}
A^{(8)}(z)=\left[  \frac{\left(  U_{0}+p\right)  \left(  U_{0}-q\right)
}{z\left(  U_{0}\right)  ^{2}}\right]  ^{1/3}=\frac{\Theta^{-1/3}}%
{\sqrt{zU_{0}}}. \label{G-}%
\end{equation}
Therefore,%
\begin{equation}
K^{(8)}(\beta,z)=\varepsilon^{1/3}\exp\left[  \varepsilon^{-1}\psi_{0}%
(z)+\ln\left(  \frac{U_{0}+p}{U_{0}-q}\right)  \beta\varepsilon^{-1/3}\right]
\mathrm{Ai}\left[  \Theta^{2/3}\beta\right]  \frac{\Theta^{-1/3}}{\sqrt
{zU_{0}}} \label{KY-}%
\end{equation}
for $y-Y^{-}(z)=O\left(  \varepsilon^{2/3}\right)  .$

\subsection{The upper part $p<z<1$ (Region IX)}

Using (\ref{beta}) in (\ref{psi+-}) we have, as $y\uparrow Y^{-}(z)$
\begin{equation}
\psi^{-}(y,z)\sim\psi_{0}(z)+\ln\left(  \frac{U_{0}+p}{U_{0}-q}\right)
\beta\varepsilon^{2/3}+\frac{2}{3}\vartheta\beta^{3/2}\varepsilon, \label{p1}%
\end{equation}%
\begin{equation}
\psi^{+}(y,z)\sim\psi_{0}(z)+\ln\left(  \frac{U_{0}+p}{U_{0}-q}\right)
\beta\varepsilon^{2/3}-\frac{2}{3}\vartheta\beta^{3/2}\varepsilon\label{p2}%
\end{equation}
and from (\ref{L+-}) we have%
\begin{equation}
L^{-}(y,z)\sim\frac{1}{\sqrt{2}}\frac{1}{\sqrt{zU_{0}\vartheta}}%
\mathrm{i}\beta^{-1/4}\varepsilon^{-1/6}, \label{l1}%
\end{equation}%
\begin{equation}
L^{+}(y,z)\sim\frac{1}{\sqrt{2}}\frac{1}{\sqrt{zU_{0}\vartheta}}\beta
^{-1/4}\varepsilon^{-1/6} \label{l2}%
\end{equation}
where%
\begin{equation}
\vartheta(z)=-\Theta(z)=\sqrt{\frac{U_{0}}{z}}\frac{1}{\left(  U_{0}+p\right)
\left(  q-U_{0}\right)  }. \label{theta}%
\end{equation}
Using (\ref{p1})-(\ref{l2}) in (\ref{left1}) we have%
\begin{gather}
K^{\left(  7\right)  }(y,z)\sim\frac{1}{2\sqrt{\pi}}\frac{1}{\sqrt
{zU_{0}\vartheta}}\beta^{-1/4}\varepsilon^{1/3}\exp\left[  \psi_{0}%
(z)\varepsilon^{-1}+\ln\left(  \frac{U_{0}+p}{U_{0}-q}\right)  \beta
\varepsilon^{-1/3}\right] \label{K30}\\
\times\left[  \frac{1}{2}\lambda^{+}(\beta,z)\exp\left(  -\frac{2}{3}%
\vartheta\beta^{3/2}\right)  +\mathrm{i}\lambda^{-}(\beta,z)\exp\left(
\frac{2}{3}\vartheta\beta^{3/2}\right)  \right] \nonumber
\end{gather}
as $y\uparrow Y^{-}(z),$ where%
\begin{equation}
\lambda^{\pm}(\beta,z)=\exp\left\{  \frac{2\pi\mathrm{i}\left[  Y^{-}\left(
z\right)  -\beta\varepsilon^{2/3}\right]  }{\varepsilon}\right\}  \pm1.
\label{lambda}%
\end{equation}

We consider the ansatz%
\begin{align}
K_{n}(x)  &  \sim K^{(9)}(\beta,z)=\varepsilon^{\nu_{9}}\exp\left[
\varepsilon^{-1}\psi_{0}(z)+\ln\left(  \frac{U_{0}+p}{U_{0}-q}\right)
\beta\varepsilon^{-1/3}\right] \label{K31}\\
&  \times\left[  \lambda^{+}(\beta,z)A^{(9)}(z)\mathrm{Ai}\left(
\vartheta^{2/3}\beta\right)  +\lambda^{-}(\beta,z)B^{(9)}(z)\mathrm{Bi}\left(
\vartheta^{2/3}\beta\right)  \right] \nonumber
\end{align}
for $y\approx Y^{-}(z),\quad p<z<1$ and unknown functions $A^{(9)}%
(z),B^{(9)}(z).$ Using (\ref{Ailarge}) in (\ref{K31}) we have, as
$\beta\rightarrow\infty$
\begin{gather}
K^{(9)}(\beta,z)\sim\varepsilon^{\nu_{9}}\exp\left[  \varepsilon^{-1}\psi
_{0}(z)+\ln\left(  \frac{U_{0}+p}{U_{0}-q}\right)  \beta\varepsilon
^{-1/3}\right] \nonumber\\
\times\left[  \lambda^{+}(\beta,z)A^{(9)}(z)\frac{1}{2\sqrt{\pi}}%
\vartheta^{-1/6}\beta^{-1/4}\exp\left(  -\frac{2}{3}\vartheta\beta
^{3/2}\right)  +\right. \label{K32}\\
\left.  \lambda^{-}(\beta,z)B^{(9)}(z)\frac{1}{\sqrt{\pi}}\vartheta
^{-1/6}\beta^{-1/4}\exp\left(  \frac{2}{3}\vartheta\beta^{3/2}\right)
\right]  .\nonumber
\end{gather}
Matching (\ref{K32}) and (\ref{K30}) we get%
\[
\nu_{9}=\frac{1}{3},\quad A^{(9)}(z)=\frac{1}{2}\frac{\vartheta^{-1/3}}%
{\sqrt{zU_{0}}},\quad B^{(9)}(z)=\frac{1}{2}\mathrm{i}\frac{\vartheta^{-1/3}%
}{\sqrt{zU_{0}}}.
\]
Thus,%
\begin{align}
K^{(9)}(\beta,z)  &  =\varepsilon^{1/3}\exp\left[  \varepsilon^{-1}\psi
_{0}(z)+\ln\left(  \frac{U_{0}+p}{U_{0}-q}\right)  \beta\varepsilon
^{-1/3}\right] \label{K33}\\
&  \times\frac{1}{2}\frac{\vartheta^{-1/3}}{\sqrt{zU_{0}}}\left[  \lambda
^{+}(\beta,z)\mathrm{Ai}\left(  \vartheta^{2/3}\beta\right)  +\mathrm{i}%
\lambda^{-}(\beta,z)\mathrm{Bi}\left(  \vartheta^{2/3}\beta\right)  \right]
.\nonumber
\end{align}

\section{The interior of $\mathbf{E}$\textbf{ }(Region X)}

We shall now find an asymptotic solution for $Y^{-}(z)<y<Y^{+}(z).$ We set
$\beta=-\widetilde{\beta},\ \widetilde{\beta}>0$ in (\ref{psi-Y-}%
)-(\ref{L-Y-}) and obtain, for $0<z<p$%
\begin{equation}
\psi^{-}(y,z)\sim\psi_{0}(z)-\ln\left(  \frac{U_{0}+p}{U_{0}-q}\right)
\widetilde{\beta}\varepsilon^{2/3}+\frac{2}{3}\Theta\mathrm{i}\widetilde
{\beta}^{3/2}\varepsilon\label{pY-1}%
\end{equation}
and%
\begin{equation}
L^{-}(y,z)\sim\frac{1}{\sqrt{2}}\frac{1}{\sqrt{zU_{0}\Theta}}e^{-\frac{1}%
{4}\pi\mathrm{i}}\widetilde{\beta}^{-1/4}\varepsilon^{-1/6}. \label{Ly-1}%
\end{equation}
Similarly, from (\ref{psi+-})-(\ref{L+-}) we obtain, for $0<z<p$%
\begin{equation}
\psi^{+}(y,z)\sim\psi_{0}(z)-\ln\left(  \frac{U_{0}+p}{U_{0}-q}\right)
\widetilde{\beta}\varepsilon^{2/3}-\frac{2}{3}\Theta\mathrm{i}\widetilde
{\beta}^{3/2}\varepsilon\label{pY-2}%
\end{equation}
and%
\begin{equation}
L^{+}(y,z)\sim\frac{1}{\sqrt{2}}\frac{1}{\sqrt{zU_{0}\Theta}}e^{\frac{1}{4}%
\pi\mathrm{i}}\widetilde{\beta}^{-1/4}\varepsilon^{-1/6}. \label{LY-2}%
\end{equation}

Introducing the function $K^{(10)}(y,z)$ defined by%
\begin{gather}
K^{(10)}(y,z)=\varepsilon^{\nu_{10}}A^{(10)}(y)\exp\left[  \varepsilon
^{-1}\psi^{-}(y,z)\right]  L^{-}(y,z)\label{K5.1}\\
+\varepsilon^{\nu_{10}}B^{(10)}(y)\exp\left[  \varepsilon^{-1}\psi
^{+}(y,z)\right]  L^{+}(y,z)\nonumber
\end{gather}
with $A^{(10)}(y),B^{(10)}(y)$ to be determined, we have%
\begin{align}
K^{(10)}(y,z)  &  \sim\varepsilon^{\nu_{10}-1/6}\frac{1}{\sqrt{2}}\frac
{1}{\sqrt{zU_{0}\Theta}}\widetilde{\beta}^{-1/4}\exp\left[  \psi_{0}%
^{-}(z)\varepsilon^{-1}-\ln\left(  \frac{U_{0}+p}{U_{0}-q}\right)
\widetilde{\beta}\varepsilon^{-1/3}\right] \label{K7.1}\\
&  \times\left[  A^{(10)}(y)\exp\left(  \frac{2}{3}\Theta\mathrm{i}%
\widetilde{\beta}^{3/2}-\frac{1}{4}\pi\mathrm{i}\right)  +B^{(10)}%
(y)\exp\left(  -\frac{2}{3}\Theta\mathrm{i}\widetilde{\beta}^{3/2}+\frac{1}%
{4}\pi\mathrm{i}\right)  \right]  .\nonumber
\end{align}
From (\ref{KY-}) we have%
\begin{align}
K^{(9)}(\beta,z)  &  \sim\varepsilon^{1/3}\frac{1}{\sqrt{\pi}}\frac{1}%
{\sqrt{zU_{0}\Theta}}\widetilde{\beta}^{-1/4}\exp\left[  \varepsilon^{-1}%
\psi_{0}^{-}(z)-\ln\left(  \frac{U_{0}+p}{U_{0}-q}\right)  \widetilde{\beta
}\varepsilon^{-1/3}\right] \label{K6.1}\\
&  \times\sin\left(  \frac{2}{3}\Theta\widetilde{\beta}^{3/2}+\frac{\pi}%
{4}\right) \nonumber
\end{align}
where we have used the asymptotic formula \cite{MR94b:00012}
\begin{equation}
\mathrm{Ai}\left(  -x\right)  \sim\frac{x^{-1/4}}{\sqrt{\pi}}\sin\left(
\frac{2}{3}x^{3/2}+\frac{\pi}{4}\right)  ,\quad x\rightarrow\infty.
\label{Ai-large}%
\end{equation}
Matching (\ref{K7.1}) and (\ref{K6.1}) we get $\nu_{10}=\frac{1}{2}$ and%
\begin{align*}
&  \frac{1}{\sqrt{2}}\left[  A^{(10)}(y)\exp\left(  \frac{2}{3}\Theta
\mathrm{i}\widetilde{\beta}^{3/2}-\frac{1}{4}\pi\mathrm{i}\right)
+B^{(10)}(y)\exp\left(  -\frac{2}{3}\Theta\mathrm{i}\widetilde{\beta}%
^{3/2}+\frac{1}{4}\pi\mathrm{i}\right)  \right] \\
&  =\frac{1}{\sqrt{\pi}}\sin\left(  \frac{2}{3}\Theta\widetilde{\beta}%
^{3/2}+\frac{\pi}{4}\right)  =\frac{1}{2\sqrt{\pi}}\left[  \exp\left(
\frac{2}{3}\Theta\mathrm{i}\widetilde{\beta}^{3/2}-\frac{1}{4}\pi
\mathrm{i}\right)  +\exp\left(  -\frac{2}{3}\Theta\mathrm{i}\widetilde{\beta
}^{3/2}+\frac{1}{4}\pi\mathrm{i}\right)  \right]
\end{align*}
from which we conclude that
\[
A^{(10)}(y)=\frac{1}{\sqrt{2\pi}}=B^{(10)}(y)
\]
and therefore%
\begin{equation}
K_{n}(x)\sim K^{(10)}(y,z)=K^{+}(y,z)+K^{-}(y,z) \label{K5}%
\end{equation}
for $Y^{-}(z)<y<Y^{+}(z).$

\subsection{Matching the interior of $\mathbf{E}$ and the upper part of the
transition layer $Y^{-}$}

We shall now verify the matching between (\ref{K5}) and (\ref{K33}). \ For
$p<z<1$ we have%
\begin{equation}
\psi^{-}(y,z)\sim(Y^{-}+\widetilde{\beta}\varepsilon^{2/3}-z)2\pi
\mathrm{i}+\psi_{0}(z)-\ln\left(  \frac{U_{0}+p}{U_{0}-q}\right)
\widetilde{\beta}\varepsilon^{2/3}-\frac{2}{3}\vartheta\mathrm{i}%
\widetilde{\beta}^{3/2}\varepsilon, \label{p5}%
\end{equation}%
\begin{equation}
L^{-}(y,z)\sim\frac{1}{\sqrt{2}}\frac{1}{\sqrt{zU_{0}\vartheta}}e^{\frac{1}%
{4}\pi\mathrm{i}}\widetilde{\beta}^{-1/4}\varepsilon^{-1/6}, \label{l5}%
\end{equation}%
\begin{equation}
\psi^{+}(y,z)\sim\psi_{0}(z)-\ln\left(  \frac{U_{0}+p}{U_{0}-q}\right)
\widetilde{\beta}\varepsilon^{2/3}+\frac{2}{3}\vartheta\mathrm{i}%
\widetilde{\beta}^{3/2}\varepsilon\label{p6}%
\end{equation}
and%
\begin{equation}
L^{+}(y,z)\sim\frac{1}{\sqrt{2}}\frac{1}{\sqrt{zU_{0}\vartheta}}e^{-\frac
{1}{4}\pi\mathrm{i}}\widetilde{\beta}^{-1/4}\varepsilon^{-1/6}. \label{l6}%
\end{equation}
Therefore,%
\begin{gather}
K^{(10)}(y,z)\sim\frac{\varepsilon^{1/3}}{2\sqrt{\pi}}\frac{\widetilde{\beta
}^{-1/4}}{\sqrt{zU_{0}\vartheta}}\exp\left[  \varepsilon^{-1}\psi_{0}%
(z)-\ln\left(  \frac{U_{0}+p}{U_{0}-q}\right)  \widetilde{\beta}%
\varepsilon^{-1/3}\right] \label{K36}\\
\times\left[  \frac{\lambda^{+}(-\widetilde{\beta},z)+\lambda^{-}%
(-\widetilde{\beta},z)}{2}\exp\left(  -\frac{2}{3}\vartheta\mathrm{i}%
\widetilde{\beta}^{3/2}+\frac{1}{4}\pi\mathrm{i}\right)  +\exp\left(  \frac
{2}{3}\vartheta\mathrm{i}\widetilde{\beta}^{3/2}-\frac{1}{4}\pi\mathrm{i}%
\right)  \right]  .\nonumber
\end{gather}
Using (\ref{Ai-large}) and \cite{MR94b:00012}%
\[
\mathrm{Bi}\left(  -x\right)  \sim\frac{x^{-1/4}}{\sqrt{\pi}}\cos\left(
\frac{2}{3}x^{3/2}+\frac{\pi}{4}\right)  ,\quad x\rightarrow\infty
\]
in (\ref{K33}) we have%
\begin{gather*}
K^{(9)}(\beta,z)\sim\varepsilon^{1/3}\exp\left[  \varepsilon^{-1}\psi
_{0}(z)-\ln\left(  \frac{U_{0}+p}{U_{0}-q}\right)  \widetilde{\beta
}\varepsilon^{-1/3}\right] \\
\times\frac{1}{2\sqrt{\pi}}\frac{\widetilde{\beta}^{-1/4}}{\sqrt
{zU_{0}\vartheta}}\left[  \lambda^{+}(-\widetilde{\beta},z)\mathrm{\sin
}\left(  \frac{2}{3}\vartheta\widetilde{\beta}^{3/2}+\frac{\pi}{4}\right)
+\mathrm{i}\lambda^{-}(-\widetilde{\beta},z)\mathrm{\cos}\left(  \frac{2}%
{3}\vartheta\widetilde{\beta}^{3/2}+\frac{\pi}{4}\right)  \right]
\end{gather*}
or%
\begin{gather}
K^{(9)}(\beta,z)\sim\varepsilon^{1/3}\exp\left[  \varepsilon^{-1}\psi
_{0}(z)-\ln\left(  \frac{U_{0}+p}{U_{0}-q}\right)  \widetilde{\beta
}\varepsilon^{-1/3}\right]  \frac{1}{2\sqrt{\pi}}\frac{\widetilde{\beta
}^{-1/4}}{\sqrt{zU_{0}\vartheta}}\label{K37}\\
\times\left[  \frac{\lambda^{+}+\lambda^{-}}{2}\exp\left(  -\frac{2}%
{3}\vartheta\mathrm{i}\widetilde{\beta}^{3/2}+\frac{1}{4}\pi\mathrm{i}\right)
+\frac{\lambda^{+}-\lambda^{-}}{2}\exp\left(  \frac{2}{3}\vartheta
\mathrm{i}\widetilde{\beta}^{3/2}-\frac{1}{4}\pi\mathrm{i}\right)  \right]
.\nonumber
\end{gather}
Since from (\ref{lambda}) we have
\[
\frac{\lambda^{+}-\lambda^{-}}{2}=1
\]
we see that (\ref{K37}) agrees with (\ref{K36}).

\subsection{Matching the interior of $\mathbf{E}$ and the corner layer at
$(p,0)$}

Using (\ref{stirling}) and the asymptotic formula \cite{MR50:2568}%
\[
H_{n}(x)\sim\sqrt{2}^{n+1}n^{n/2}\exp\left(  \frac{x^{2}}{2}-\frac{n}%
{2}\right)  \cos\left(  \sqrt{2n}x-\frac{n\pi}{2}\right)  ,\quad
n\rightarrow\infty
\]
in (\ref{Hn}) we have, as $n\rightarrow\infty$%
\begin{equation}
K_{n}^{\left(  2\right)  }(\eta)\sim\frac{1}{\sqrt{n\pi}}\exp\left\{
\frac{\eta^{2}}{2}+\frac{n}{2}\left[  1+\ln\left(  \frac{pq}{\varepsilon
n}\right)  \right]  \right\}  \cos\left(  \sqrt{2n}\eta-\frac{n\pi}{2}\right)
. \label{K2large}%
\end{equation}

Using (\ref{eta}) in (\ref{psi+-}), with $z=n\varepsilon$ we get%
\begin{equation}
\frac{1}{\varepsilon}\psi^{\pm}\sim\frac{n}{2}\left[  1+\ln\left(  \frac
{pq}{\varepsilon n}\right)  \pm\pi\mathrm{i}\right]  \mp\sqrt{2n}%
\eta\mathrm{i}+\frac{\eta^{2}}{2}.\label{12}%
\end{equation}
Similarly from (\ref{L+-}) we find%
\begin{equation}
L^{+}\sim\frac{1}{\sqrt{2n}},\quad L^{-}\sim\frac{1}{\sqrt{2n}}.\label{14}%
\end{equation}
Using (\ref{12})-(\ref{14}) in (\ref{K5}) we obtain%
\[
K^{(10)}\sim\frac{1}{\sqrt{\pi n}}\exp\left\{  \frac{\eta^{2}}{2}+\frac{n}%
{2}\left[  1+\ln\left(  \frac{pq}{\varepsilon n}\right)  \right]  \right\}
\cos\left(  \sqrt{2n}\eta-\frac{n\pi}{2}\right)
\]
which is in agreement with (\ref{K2large}).

\section{The boundary layer at $z=1$ (Region XI)}

We now consider solutions of (\ref{recurrence}) with $n\approx N.$ We
introduce the variable $j$ and the function $\Upsilon_{j}^{(11)}(y)$ defined
by
\begin{equation}
n=N-j,\quad j\in\mathbb{Z},\quad j\geq-1, \label{j}%
\end{equation}%
\begin{equation}
K_{n}(x)=\binom{N}{N-n}\left(  -p\right)  ^{n}\Upsilon_{N-n}^{(11)}\left(
\frac{x}{N}\right)  . \label{z=1(1)}%
\end{equation}
Using (\ref{z=1(1)}) in (\ref{recurrence}) we have, to leading order,%
\[
q\Upsilon_{j+1}^{(11)}-\left(  q-y\right)  \Upsilon_{j}^{(11)}=0
\]
which can be solved to obtain%
\begin{equation}
\Upsilon_{j}^{(11)}(y)=A^{(11)}\left(  y\right)  \left(  1-\frac{y}{q}\right)
^{j} \label{z=1(2)}%
\end{equation}
where $A^{(11)}\left(  y\right)  $ is a function to be determined. Thus, we
have%
\begin{equation}
K_{n}(x)\sim A^{(11)}\left(  y\right)  \binom{N}{j}\left(  -p\right)
^{N-j}\left(  1-\frac{y}{q}\right)  ^{j}. \label{z=1(3)}%
\end{equation}

We now introduce the function $R_{j}^{(11)}(y)$ defined by%
\begin{equation}
K_{n}(x)=\binom{x}{n}R_{N-n}^{(11)}\left(  \frac{x}{N}\right)  \label{z=1(4)}%
\end{equation}
and note that from (\ref{KN}) we have
\begin{equation}
R_{-1}^{(11)}(y)=1. \label{R60}%
\end{equation}
Using (\ref{j}) in (\ref{recurrence}) we obtain, to leading order,%
\[
(q-y)R_{j+1}^{(11)}-(1-y)R_{j}^{(11)}=0
\]
which together with (\ref{R60}) implies%
\begin{equation}
R_{j}^{(11)}(y)=\left(  \frac{1-y}{q-y}\right)  ^{j+1},\quad y\neq q.
\label{R81}%
\end{equation}
Therefore,
\begin{equation}
K_{n}(x)\sim\binom{Ny}{n}\left(  \frac{1-y}{q-y}\right)  ^{j+1},\quad
y\not \approx q. \label{R82}%
\end{equation}
From (\ref{z=1(4)}) and (\ref{R82}) we get, for $y\not \approx q,$
\begin{equation}
K_{n}(x)\sim K_{j}^{(11)}(y)=A^{(11)}\left(  y\right)  \binom{N}{j}\left(
-p\right)  ^{N-j}\left(  1-\frac{y}{q}\right)  ^{j}+\binom{Ny}{N-j}\left(
\frac{1-y}{q-y}\right)  ^{j+1}. \label{K6}%
\end{equation}

\subsection{Matching with the upper corners}

We shall now determine the function $A^{(11)}\left(  y\right)  $ by matching
(\ref{K6}) with (\ref{leftcorner}). As $N\rightarrow\infty$ we get%
\begin{gather}
K_{j}^{(11)}(y)\sim A^{(11)}\left(  y\right)  \frac{1}{\sqrt{2\pi j}}%
\exp\left\{  \varepsilon^{-1}\ln\left(  -p\right)  +\left[  1-\ln\left(
j\varepsilon\right)  -\ln(-p)+\ln\left(  1-\frac{y}{q}\right)  \right]
j\right\} \label{K8Large}\\
+\varepsilon^{1/2}\sqrt{\frac{2}{\pi}}\frac{\sqrt{\left(  1-y\right)  y}}%
{y-q}\sin(N\pi y)\exp\left\{  \frac{y\ln(y)+(1-y)\ln\left(  1-y\right)
+\pi\mathrm{i}}{\varepsilon}-j\ln\left(  y-q\right)  \right\}  .\nonumber
\end{gather}

Using (\ref{j}) in (\ref{psi+-}) we have as $\varepsilon\rightarrow0,$
\ $0<y<q$%
\begin{equation}
\psi^{+}(y,z)\sim\ln(-p)+y\ln\left(  \frac{q}{p}\right)  -y\pi\mathrm{i}%
+\left[  1-\ln\left(  j\varepsilon\right)  -\ln(-p)+\ln\left(  1-\frac{y}%
{q}\right)  \right]  j\varepsilon\label{psipz=1}%
\end{equation}%
\begin{equation}
\psi^{-}(y,z)\sim y\ln(y)+(1-y)\ln\left(  1-y\right)  +(1-y)\pi\mathrm{i}%
-j\ln(y-q)\varepsilon\label{psimz=1}%
\end{equation}
and from (\ref{L+-}) we get%
\begin{equation}
L^{+}(y,z)\sim\frac{1}{\sqrt{j\varepsilon}},\quad L^{-}(y,z)\sim\sqrt
{y(1-y)}\frac{1}{q-y}\mathrm{i}. \label{Lz=1}%
\end{equation}
Using (\ref{psipz=1})-(\ref{Lz=1}) in (\ref{leftcorner}) we have%
\begin{gather}
K^{\left(  7\right)  }(y,z)\sim\cos\left(  \frac{\pi y}{\varepsilon}\right)
\frac{1}{\sqrt{2\pi j}}\exp\left\{  \frac{\ln(-p)+y\ln\left(  \frac{q}%
{p}\right)  }{\varepsilon}+\left[  1-\ln\left(  j\varepsilon\right)
-\ln(-p)+\ln\left(  1-\frac{y}{q}\right)  \right]  j\right\}
\label{leftcornerz=1}\\
+\sqrt{y(1-y)}\frac{1}{y-q}\sin\left(  \frac{\pi y}{\varepsilon}\right)
\sqrt{\frac{2\varepsilon}{\pi}}\exp\left[  \frac{y\ln(y)+(1-y)\ln\left(
1-y\right)  +\pi\mathrm{i}}{\varepsilon}-j\ln(y-q)\right]  .\nonumber
\end{gather}

Matching (\ref{K8Large}) and \ref{leftcornerz=1} we obtain%
\begin{equation}
A^{(11)}(y)=\exp\left[  y\ln\left(  \frac{q}{p}\right)  \varepsilon
^{-1}\right]  \cos\left(  \frac{\pi y}{\varepsilon}\right)  . \label{a(y)}%
\end{equation}
Hence, we conclude that%
\begin{equation}
K_{j}^{(11)}(y)=\binom{N}{j}\left(  -p\right)  ^{N-j}\left(  \frac{q}%
{p}\right)  ^{Ny}\cos\left(  N\pi y\right)  \left(  1-\frac{y}{q}\right)
^{j}+\binom{Ny}{N-j}\left(  \frac{1-y}{q-y}\right)  ^{j+1} \label{upper}%
\end{equation}
for $z=1-O(\varepsilon)$ and $y\not \approx q.$

\subsection{The corner layer at $(q,1)$ (Region XII)}

The approximation (\ref{upper}) ceases to be valid for $y\approx q.$
Therefore, we need to find another expression, which holds in a neighborhood
of the point $(q,1).$ We introduce the new variable $\xi$ defined by%
\begin{equation}
y=q+\xi\sqrt{2pq\varepsilon},\quad\xi=O(1). \label{xsi}%
\end{equation}
Using (\ref{xsi}) in (\ref{upper}) we have, as $\varepsilon\rightarrow0$%
\begin{gather}
K_{j}^{(11)}(y)\sim\exp\left[  \frac{p\ln(p)+q\ln(q)}{\varepsilon}+\xi
\sqrt{\frac{2pq}{\varepsilon}}\ln\left(  \frac{q}{p}\right)  -\frac{j}{2}%
\ln(pq\varepsilon)\right] \label{z=1,y=q1}\\
\times\left[  \frac{\left(  \sqrt{2}\xi\right)  ^{j}}{j!}\cos\left(
\frac{p\pi}{\varepsilon}-\xi\pi\sqrt{\frac{2pq}{\varepsilon}}\right)
-\sqrt{\frac{2}{\pi}}e^{\xi^{2}}\left(  \sqrt{2}\xi\right)  ^{-j-1}\sin\left(
\frac{p\pi}{\varepsilon}-\xi\pi\sqrt{\frac{2pq}{\varepsilon}}\right)  \right]
.\nonumber
\end{gather}
Equation (\ref{z=1,y=q1}) suggests that we define a new function $R_{j}%
^{(12)}(\xi)$ by%
\begin{equation}
K_{n}(x)=\exp\left[  \frac{p\ln(p)+q\ln(q)}{\varepsilon}+\xi\sqrt{\frac
{2pq}{\varepsilon}}\ln\left(  \frac{q}{p}\right)  -\frac{j}{2}\ln
(pq\varepsilon)\right]  R_{j}^{(12)}(\xi). \label{z=1y=q2}%
\end{equation}
Using (\ref{z=1y=q2}) in (\ref{recurrence}) we obtain, to leading order,%
\[
\left(  j+1\right)  R_{j+1}^{(12)}-\sqrt{2}\xi R_{j}+R_{j-1}^{(12)}=0
\]
which has the independent solutions \cite{MR97k:01072}%
\[
R_{j}^{(12)}\left(  \xi\right)  =\frac{1}{j!}\mathrm{D}_{j}\left(  \sqrt{2}%
\xi\right)  ,\text{ \ and \ }R_{j}^{(12)}\left(  \xi\right)  =\mathrm{D}%
_{-j-1}\left(  \pm\sqrt{2}\mathrm{i}\xi\right)  \left(  \pm\mathrm{i}\right)
^{j}.
\]
Using (\ref{Dularge}) and matching with (\ref{z=1,y=q1}) we get%
\begin{gather}
K_{n}(x)\sim K_{j}^{\left(  12\right)  }(\xi)=\exp\left[  \frac{p\ln
(p)+q\ln(q)}{\varepsilon}+\xi\sqrt{\frac{2pq}{\varepsilon}}\ln\left(  \frac
{q}{p}\right)  -\frac{j}{2}\ln(pq\varepsilon)\right]  \exp\left(  \frac
{\xi^{2}}{2}\right) \label{z=1y=q3}\\
\times\left[  \frac{1}{j!}\mathrm{D}_{j}\left(  \sqrt{2}\xi\right)
\cos\left(  \frac{p\pi}{\varepsilon}-\xi\pi\sqrt{\frac{2pq}{\varepsilon}%
}\right)  -\frac{1}{\sqrt{2\pi}}\Lambda_{j}\left(  \xi\right)  \sin\left(
\frac{p\pi}{\varepsilon}-\xi\pi\sqrt{\frac{2pq}{\varepsilon}}\right)  \right]
,\nonumber
\end{gather}
where the function $\Lambda_{j}:\mathbb{R}\rightarrow\mathbb{R}$ is defined
by
\begin{equation}
\Lambda_{j}\left(  \xi\right)  =\mathrm{i}^{j+1}\left[  \mathrm{D}%
_{-j-1}\left(  \sqrt{2}\mathrm{i}\xi\right)  +\left(  -1\right)
^{j+1}\mathrm{D}_{-j-1}\left(  -\sqrt{2}\mathrm{i}\xi\right)  \right]  .
\label{Lambda}%
\end{equation}

\subsubsection{Matching the corner at $(q,1)$ and the interior of $\mathbf{E}%
$}

Finally, we verify the matching between (\ref{K5}) and (\ref{z=1y=q3}). Using
(\ref{j}) and (\ref{xsi}) in (\ref{psi+-}) we have,%
\begin{gather*}
\frac{\psi^{+}}{\varepsilon}\sim\frac{q\ln(q)+p\ln(p)+p\pi\mathrm{i}%
}{\varepsilon}+\xi\sqrt{2pq}\left[  \ln\left(  \frac{q}{p}\right)
-\pi\mathrm{i}\right]  \varepsilon^{-1/2}\\
+\frac{j}{2}\left[  1-\ln(pqj\varepsilon)-\pi\mathrm{i}\right]  +\sqrt{2j}%
\xi\mathrm{i}+\frac{\xi^{2}}{2}%
\end{gather*}%
\begin{gather*}
\frac{\psi^{-}}{\varepsilon}\sim\frac{q\ln(q)+p\ln(p)-p\pi\mathrm{i}%
}{\varepsilon}+\xi\sqrt{2pq}\left[  \ln\left(  \frac{q}{p}\right)
+\pi\mathrm{i}\right]  \varepsilon^{-1/2}\\
+\frac{j}{2}\left[  1-\ln(pqj\varepsilon)+\pi\mathrm{i}\right]  -\sqrt{2j}%
\xi\mathrm{i}+\frac{\xi^{2}}{2}.
\end{gather*}
Similarly, from (\ref{L+-}) we get%
\[
L^{+}\sim\frac{1}{\sqrt{2j\varepsilon}},\quad L^{-}\sim\frac{1}{\sqrt
{2j\varepsilon}}%
\]
and therefore%
\begin{gather}
K^{\left(  10\right)  }(y,z)\sim\frac{1}{\sqrt{\pi j}}\exp\left[  \frac
{q\ln(q)+p\ln(p)}{\varepsilon}+\xi\sqrt{2pq}\ln\left(  \frac{q}{p}\right)
\varepsilon^{-1/2}+\frac{\xi^{2}}{2}\right] \label{z=1,y=q4}\\
\times\exp\left\{  \frac{j}{2}\left[  1-\ln(pqj\varepsilon)\right]  \right\}
\cos\left[  \frac{p\pi}{\varepsilon}-\xi\sqrt{2pq}\pi\varepsilon^{-1/2}%
+\sqrt{2j}\xi-\frac{j\pi}{2}\right]  .\nonumber
\end{gather}
Using the formula \cite{MR97k:01072}
\[
\mathrm{D}_{n}\left(  \xi\right)  \sim\sqrt{2}\exp\left\{  \frac{n}{2}\left[
\ln(n)-1\right]  \right\}  \cos\left[  \sqrt{n}\xi-\frac{n\pi}{2}\right]
,\quad n\rightarrow\infty
\]
in (\ref{z=1y=q3}) yields%
\begin{gather}
K_{j}^{\left(  12\right)  }(\xi)\sim\exp\left[  \frac{p\ln(p)+q\ln
(q)}{\varepsilon}+\xi\sqrt{\frac{2pq}{\varepsilon}}\ln\left(  \frac{q}%
{p}\right)  -\frac{j}{2}\ln(pq\varepsilon)\right]  \frac{1}{\sqrt{\pi j}}%
\exp\left\{  \frac{\xi^{2}}{2}+\frac{j}{2}\left[  1-\ln(j)\right]  \right\}
\label{z=1y=q5}\\
\times\left[  \cos\left(  \sqrt{j}\xi-\frac{j\pi}{2}\right)  \cos\left(
\frac{p\pi}{\varepsilon}-\xi\pi\sqrt{\frac{2pq}{\varepsilon}}\right)
-\sin\left(  \sqrt{j}\xi-\frac{j\pi}{2}\right)  \sin\left(  \frac{p\pi
}{\varepsilon}-\xi\pi\sqrt{\frac{2pq}{\varepsilon}}\right)  \right]
,\nonumber
\end{gather}
where we have used
\[
\Lambda_{j}\left(  \xi\right)  \sim\sqrt{\frac{2}{j}}\exp\left\{  \frac{j}%
{2}\left[  1-\ln(j)\right]  \right\}  \sin\left(  \sqrt{j}\xi-\frac{j\pi}%
{2}\right)  .
\]
Equations (\ref{z=1,y=q4}) and (\ref{z=1y=q5}) are identical, after regrouping terms.

\section{Summary and numerical results}

\begin{figure}[ptb]
\begin{center}
\rotatebox{270} {\resizebox{12cm}{!}{\includegraphics{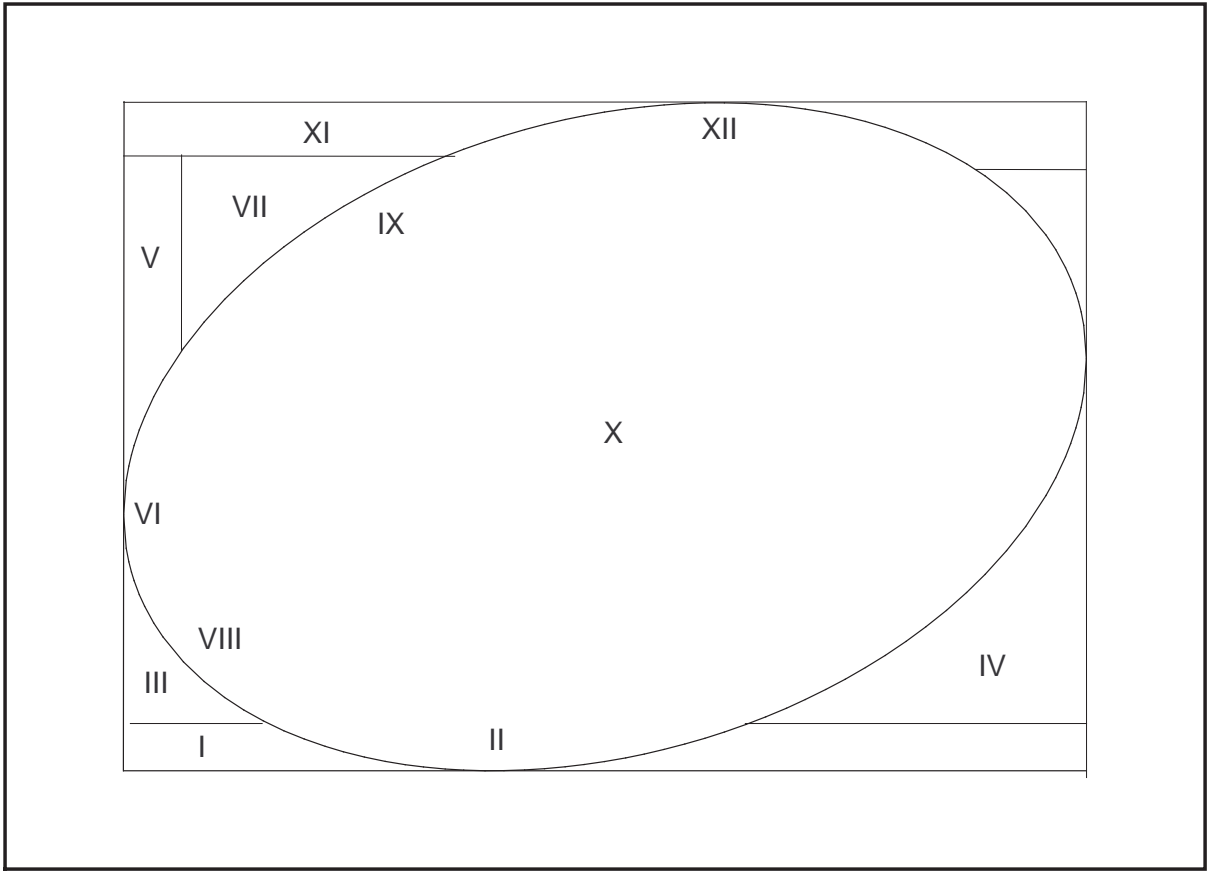}}}
\end{center}
\caption{A sketch of the different asymptotic regions.}%
\label{regions}%
\end{figure}

Below we summarize our results for the various asymptotic approximations to
$K_{n}(x)$ as $N\rightarrow\infty,$ with
\[
\varepsilon=\frac{1}{N},\quad x=\frac{y}{\varepsilon},\quad n=\frac
{z}{\varepsilon}\quad0\leq y,z\leq1.
\]
(see Figure \ref{regions}).

\begin{enumerate}
\item Region I: $\ n=O(1),$ $0\leq y\leq1,$ $y\not \approx p.$%
\[
K_{n}(x)\sim K_{n}^{(1)}(y)=\frac{\varepsilon^{-n}}{n!}\left(  y-p\right)
^{n}.
\]

(see Figure \ref{RegionI}).

\begin{figure}[ptb]
\begin{center}
\rotatebox{270} {\resizebox{12cm}{!}{\includegraphics{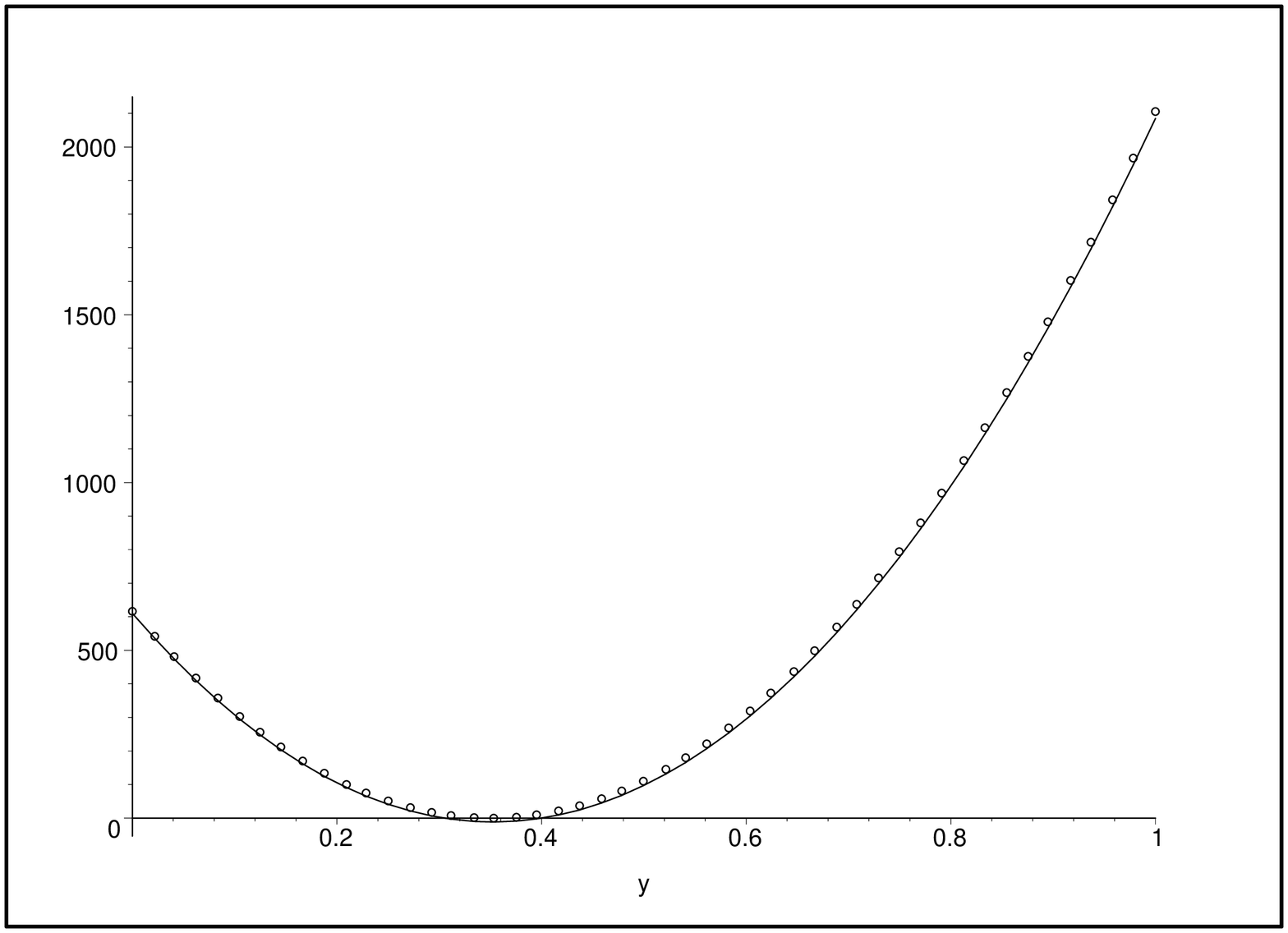}}}
\end{center}
\caption{A comparison of $K_{n}(x)$ (solid curve) and $K_{n}^{(1)}(y)$ (ooo)
for $n=2$ with $\varepsilon=0.01$ and $q=0.64894783$.}%
\label{RegionI}%
\end{figure}

\item Region II: \ $n=O(1),$ $y\approx p,$ $y=p+\eta\sqrt{2pq\varepsilon},$
$\eta=O(1).$
\[
K_{n}(x)\sim K_{n}^{\left(  2\right)  }(\eta)=\frac{\varepsilon^{-n/2}}%
{n!}\left(  \frac{pq}{2}\right)  ^{n/2}H_{n}\left(  \eta\right)  ,
\]
where $H_{n}\left(  \eta\right)  $ is the Hermite polynomial (see Figure
\ref{RegionII}).

\begin{figure}[ptb]
\begin{center}
\rotatebox{270} {\resizebox{12cm}{!}{\includegraphics{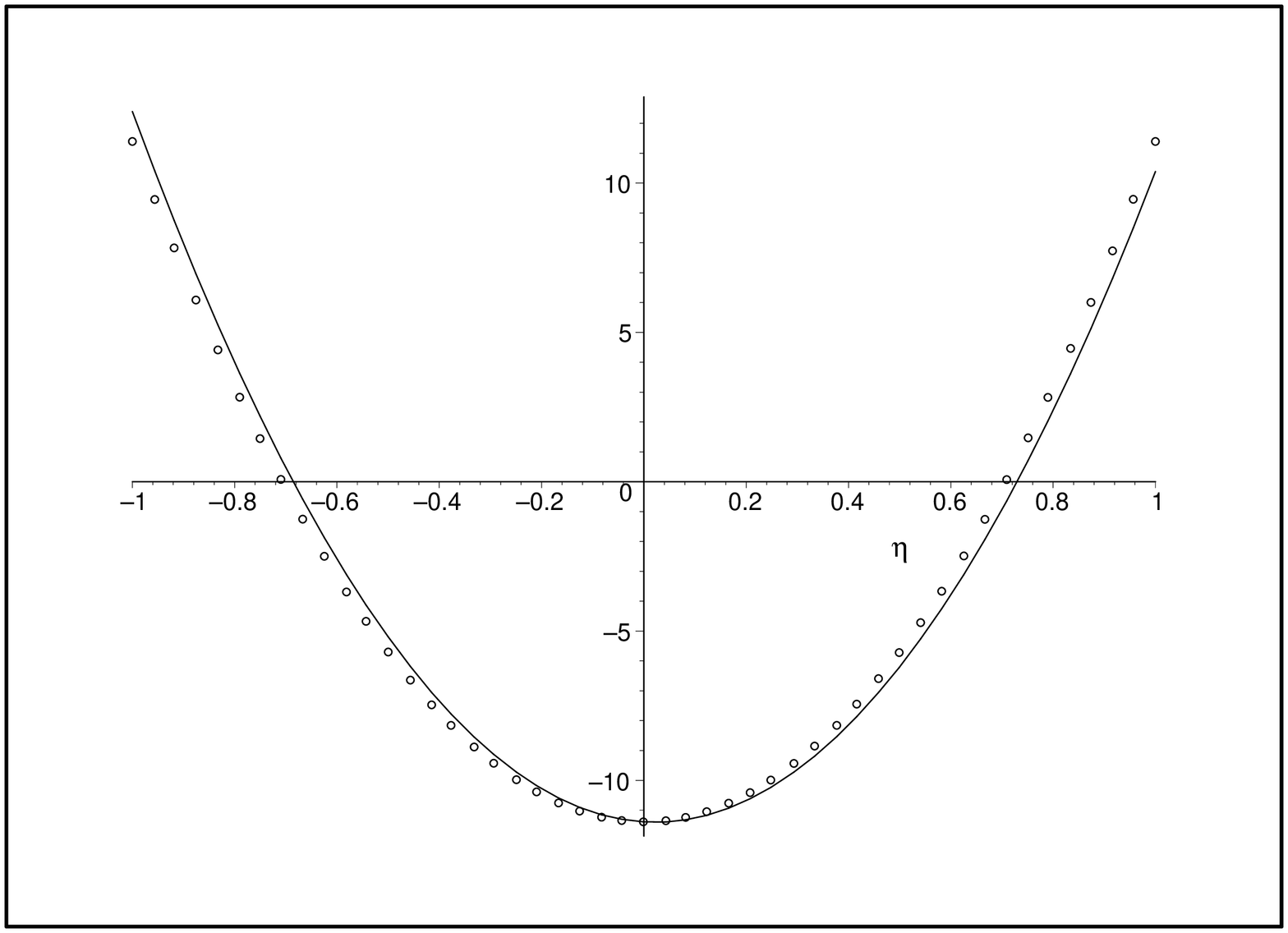}}}
\end{center}
\caption{A comparison of $K_{n}(x)$ (solid curve) and $K_{n}^{(2)}(\eta)$
(ooo) for $n=2$ with $\varepsilon=0.01$ and $q=0.64894783$.}%
\label{RegionII}%
\end{figure}

\item Region III: \ $0\leq y<Y^{-}(z),$ $0<z<p,$ where
\[
Y^{\pm}(z)=p+\left(  q-p\right)  z\pm2zU_{0},\quad U_{0}(z)=\sqrt
{\frac{pq(1-z)}{z}}.
\]%
\[
K_{n}(x)\sim K^{-}(y,z)=\varepsilon^{1/2}\frac{1}{\sqrt{2\pi}}\exp\left[
\varepsilon^{-1}\psi^{-}(y,z)\right]  L^{-}(y,z),
\]
with%
\[
\psi^{\pm}(y,z)=\ln\left[  \left(  U^{\pm}\right)  ^{z-1}\left(  U^{\pm
}-p\right)  ^{1-y}\left(  U^{\pm}+q\right)  ^{y}\right]  ,\quad L^{\pm
}(y,z)=\sqrt{\frac{(U^{\pm}-p)(U^{\pm}+q)}{z\left[  \left(  U^{\pm}\right)
^{2}-U_{0}^{2}\right]  }}%
\]
and%
\[
U^{\pm}(y,z)=-\frac{1}{2}\left(  \frac{p-y}{z}+q-p\right)  \pm\frac{1}{2}%
\sqrt{\left(  \frac{p-y}{z}+q-p\right)  ^{2}-4\left(  U_{0}\right)  ^{2}}
\]
(see Figure \ref{RegionIII}).

\begin{figure}[ptb]
\begin{center}
\rotatebox{270} {\resizebox{12cm}{!}{\includegraphics{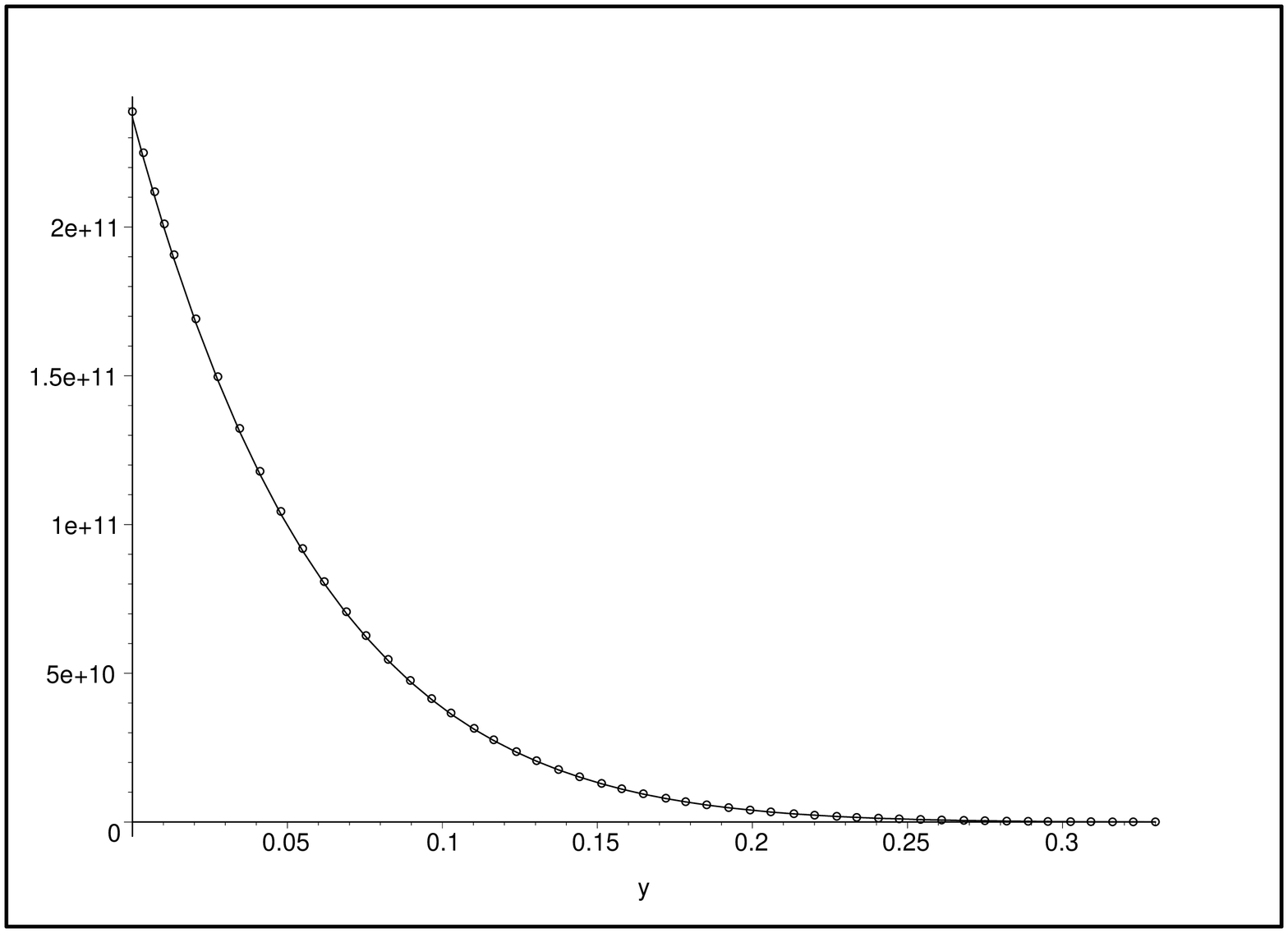}}}
\end{center}
\caption{A comparison of $K_{n}(x)$ (solid curve) and $K^{-}(y,z)$ (ooo) for
$n=10$ with $\varepsilon=0.01$ and $q=0.34894783$.}%
\label{RegionIII}%
\end{figure}

\item Region IV: $Y^{+}(z)<y\leq1,$ $0<z<q.$%
\[
K_{n}(x)\sim K^{+}(y,z)=\varepsilon^{1/2}\frac{1}{\sqrt{2\pi}}\exp\left[
\varepsilon^{-1}\psi^{+}(y,z)\right]  L^{+}(y,z).
\]

(see Figure \ref{RegionIV}).

\begin{figure}[ptb]
\begin{center}
\rotatebox{270} {\resizebox{12cm}{!}{\includegraphics{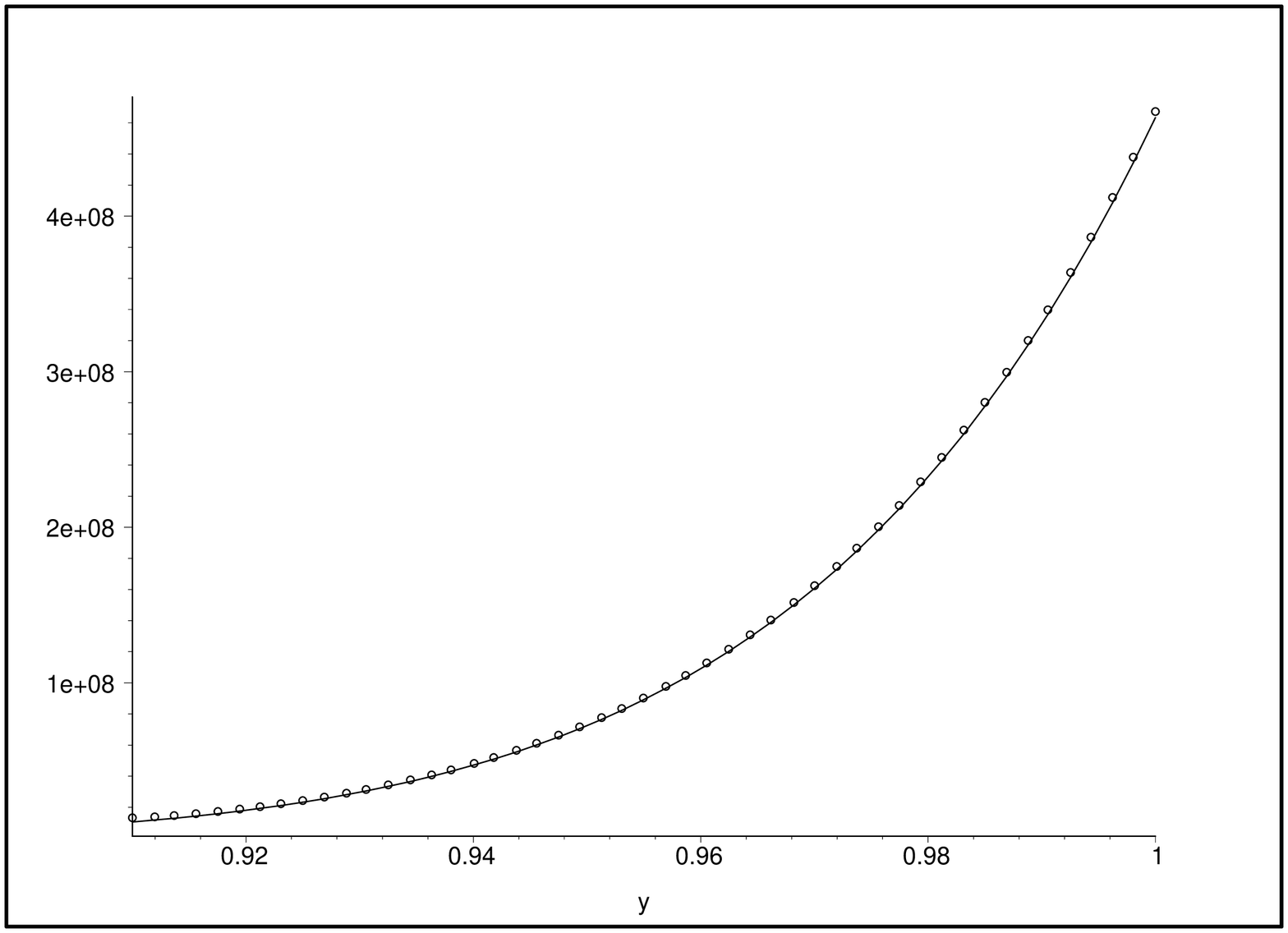}}}
\end{center}
\caption{A comparison of $K_{n}(x)$ (solid curve) and $K^{+}(y,z)$ (ooo) for
$n=10$ with $\varepsilon=0.01$ and $q=0.34894783$.}%
\label{RegionIV}%
\end{figure}

\item Region V: $\ x=O(1),$ $p<z<1.$%
\begin{gather*}
K_{n}(x)\sim K^{(5)}\left(  x,z\right)  =\frac{\varepsilon^{1/2}}{\sqrt{2\pi
}\sqrt{z\left(  1-z\right)  }}\cos(\pi x)\left(  \frac{z-p}{p}\right)
^{x}\exp\left[  \frac{\phi_{0}(z)}{\varepsilon}\right] \\
-\frac{\varepsilon}{\pi}\frac{x}{z-p}\Gamma(x)\sin(\pi x)\left(
\frac{q\varepsilon}{z-p}\right)  ^{x}\exp\left[  \frac{(z-1)\ln(q)+\pi
\mathrm{i}z}{\varepsilon}\right]
\end{gather*}
where%
\[
\phi_{0}(z)=(z-1)\ln(1-z)-z\ln(z)+z\ln(-p)
\]
and $\Gamma(x)$ is the Gamma function (see Figure \ref{RegionV}).

\begin{figure}[ptb]
\begin{center}
\rotatebox{270} {\resizebox{12cm}{!}{\includegraphics{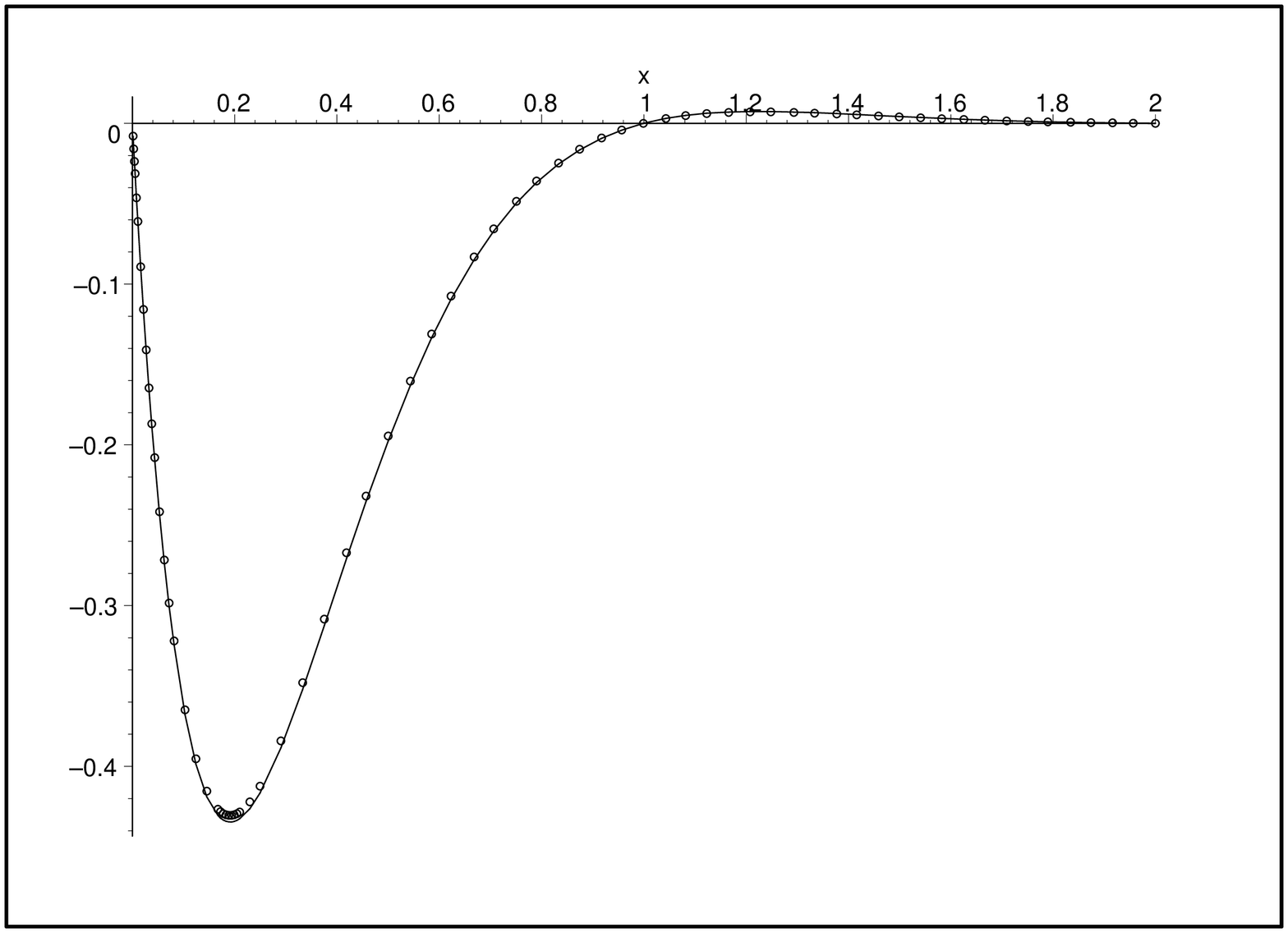}}}
\end{center}
\caption{A comparison of $K_{n}(x)$ (solid curve) and $K^{(5)}\left(
x,z\right)  $ (ooo) for $n=80$ with $\varepsilon=0.01$ and $q=0.74894783$.}%
\label{RegionV}%
\end{figure}.

\item Region VI: $\ x=O(1),$ $z\approx p,$ $z=p-u\sqrt{pq\varepsilon},$
$u=O(1).$%
\begin{align*}
K_{n}(x)  &  \sim K^{(6)}(x,u)=\frac{\varepsilon^{1/2}}{\sqrt{2\pi pq}}\left[
\sqrt{\frac{q\varepsilon}{p}}\right]  ^{x}\mathrm{D}_{x}(u)\\
&  \times\exp\left[  \frac{\pi\mathrm{i}p-q\ln\left(  q\right)  }{\varepsilon
}+\frac{u\sqrt{pq}\pi\mathrm{i-}u\sqrt{pq}\ln\left(  q\right)  }%
{\sqrt{\varepsilon}}-\frac{u^{2}}{4}\right]  ,
\end{align*}
where $\mathrm{D}_{x}(u)$ is the parabolic cylinder function (see Figure
\ref{RegionVI}).

\begin{figure}[ptb]
\begin{center}
\rotatebox{270} {\resizebox{12cm}{!}{\includegraphics{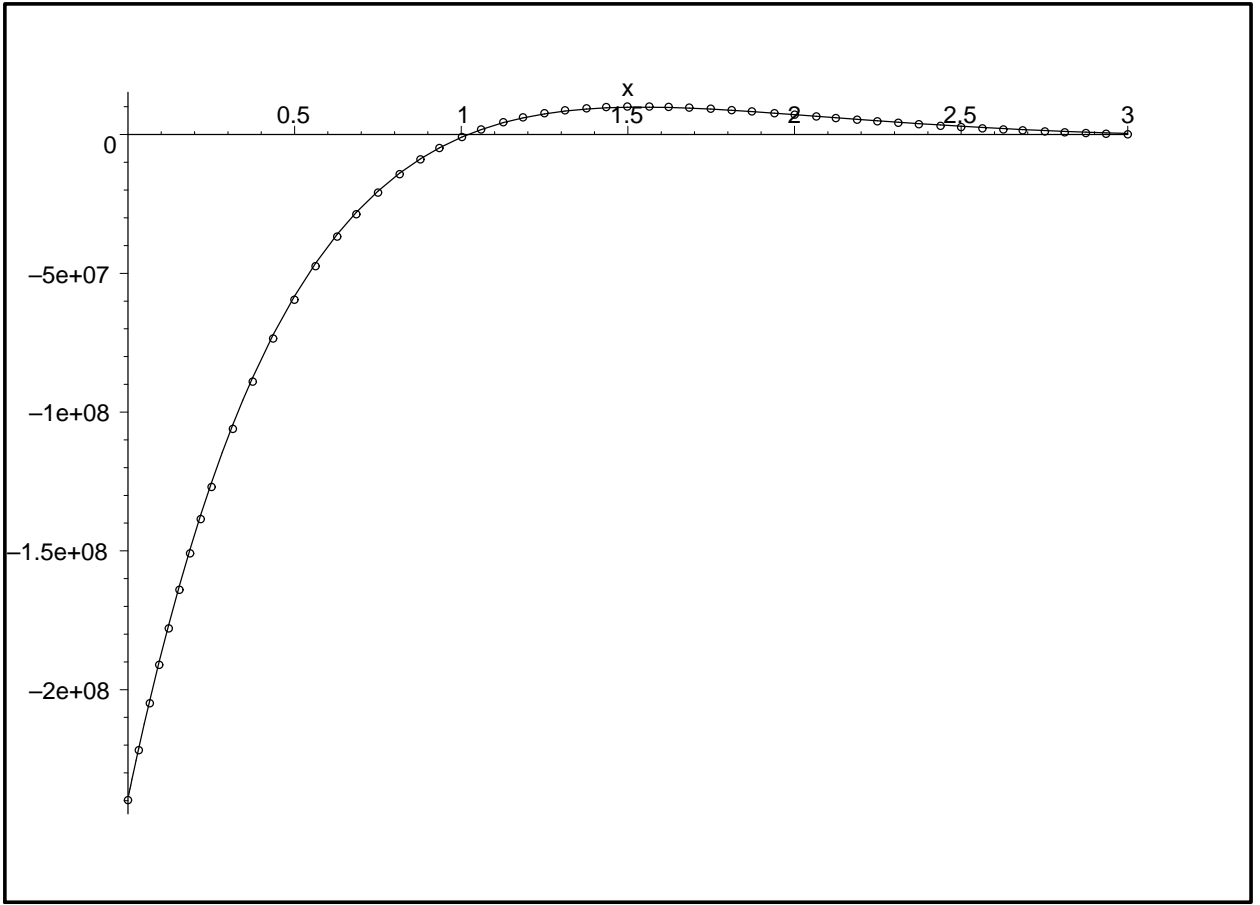}}}
\end{center}
\caption{A comparison of $K_{n}(x)$ (solid curve) and $K^{(6)}(x,u)$ (ooo) for
$n=25$ with $\varepsilon=0.01$ and $q=0.74894783$. Note that with this choice
of parameters, $u=0.024265$.}%
\label{RegionVI}%
\end{figure}

\item Region VII: $\ 0\ll y<Y^{-}(z),$ $p<z<1.$%
\[
K_{n}(x)\sim K^{\left(  7\right)  }(y,z)=\exp\left(  \frac{\pi\mathrm{i}%
y}{\varepsilon}\right)  \left[  \cos\left(  \frac{\pi y}{\varepsilon}\right)
K^{+}(y,z)+2\mathrm{i}\sin\left(  \frac{\pi y}{\varepsilon}\right)
K^{-}(y,z)\right]
\]
(see Figure \ref{RegionVII}).

\begin{figure}[ptb]
\begin{center}
\rotatebox{270} {\resizebox{12cm}{!}{\includegraphics{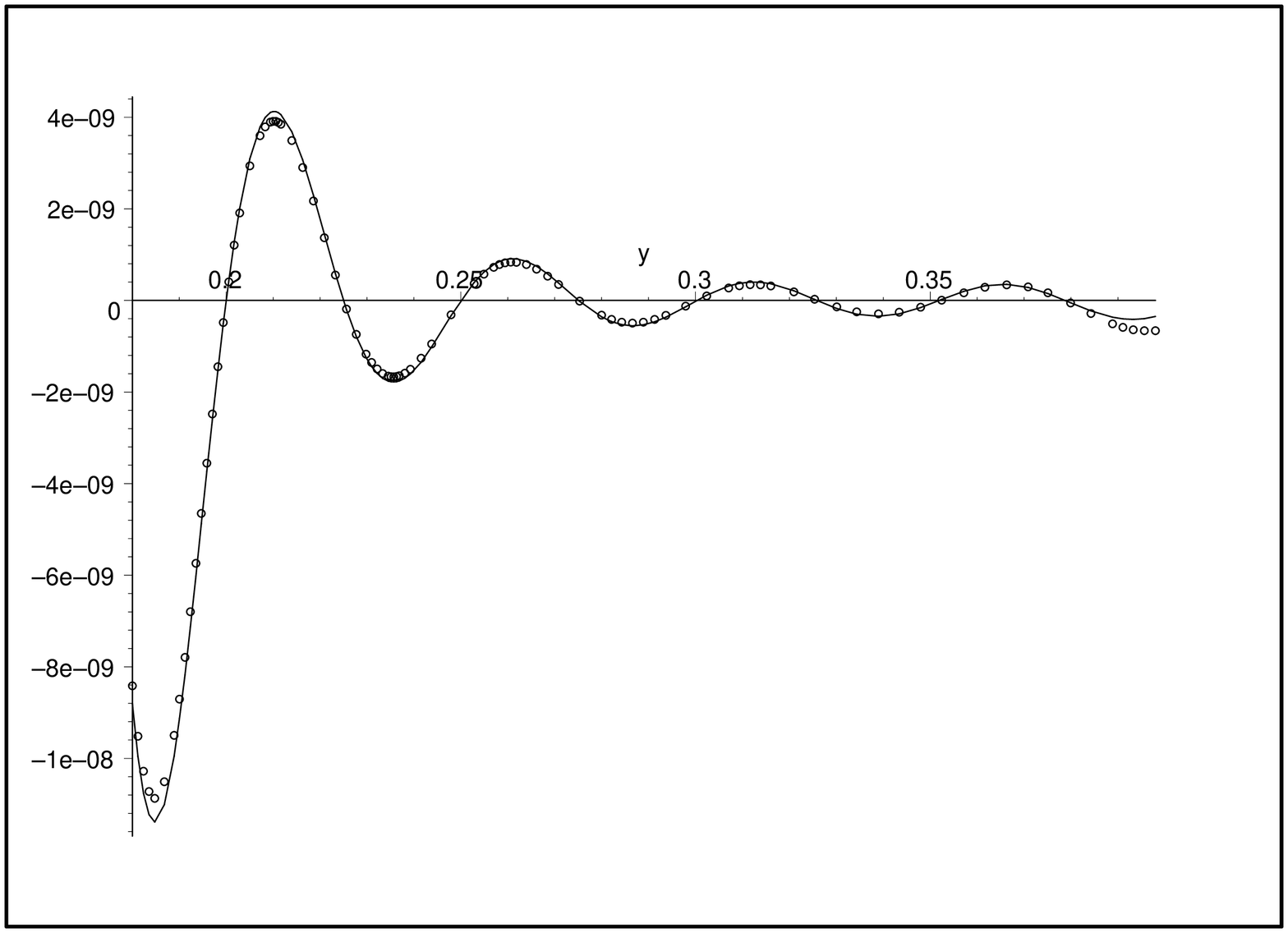}}}
\end{center}
\caption{A comparison of $K_{n}(x)$ (solid curve) and $K^{(7)}(y,z)$ (ooo) for
$n=35$ with $\varepsilon=0.025$ and $q=0.74894783$.}%
\label{RegionVII}%
\end{figure}

\item Region VIII: \ $y\approx Y^{-}(z),$ $0<z<p,$ $y=Y^{-}(z)-\beta
\varepsilon^{2/3},$ $\beta=O(1)$.%
\[
K_{n}(x)\sim K^{(8)}(\beta,z)=\varepsilon^{1/3}\exp\left[  \varepsilon
^{-1}\psi_{0}(z)+\ln\left(  \frac{U_{0}+p}{U_{0}-q}\right)  \beta
\varepsilon^{-1/3}\right]  \mathrm{Ai}\left[  \Theta^{2/3}\beta\right]
\frac{\Theta^{-1/3}}{\sqrt{zU_{0}}},
\]
where%
\[
\psi_{0}(z)=z\pi\mathrm{i}+(z-1)\ln\left(  U_{0}\right)  +Y^{-}(z)\ln\left(
U_{0}-q\right)  +\left[  1-Y^{-}(z)\right]  \ln\left(  U_{0}+p\right)  ,
\]%
\[
\Theta(z)=\sqrt{\frac{U_{0}}{z}}\frac{1}{\left(  U_{0}+p\right)  \left(
U_{0}-q\right)  }%
\]
and $\mathrm{Ai}\left(  \cdot\right)  $ is the Airy function (see Figure
\ref{RegionVIII}).

\begin{figure}[ptb]
\begin{center}
\rotatebox{270} {\resizebox{12cm}{!}{\includegraphics{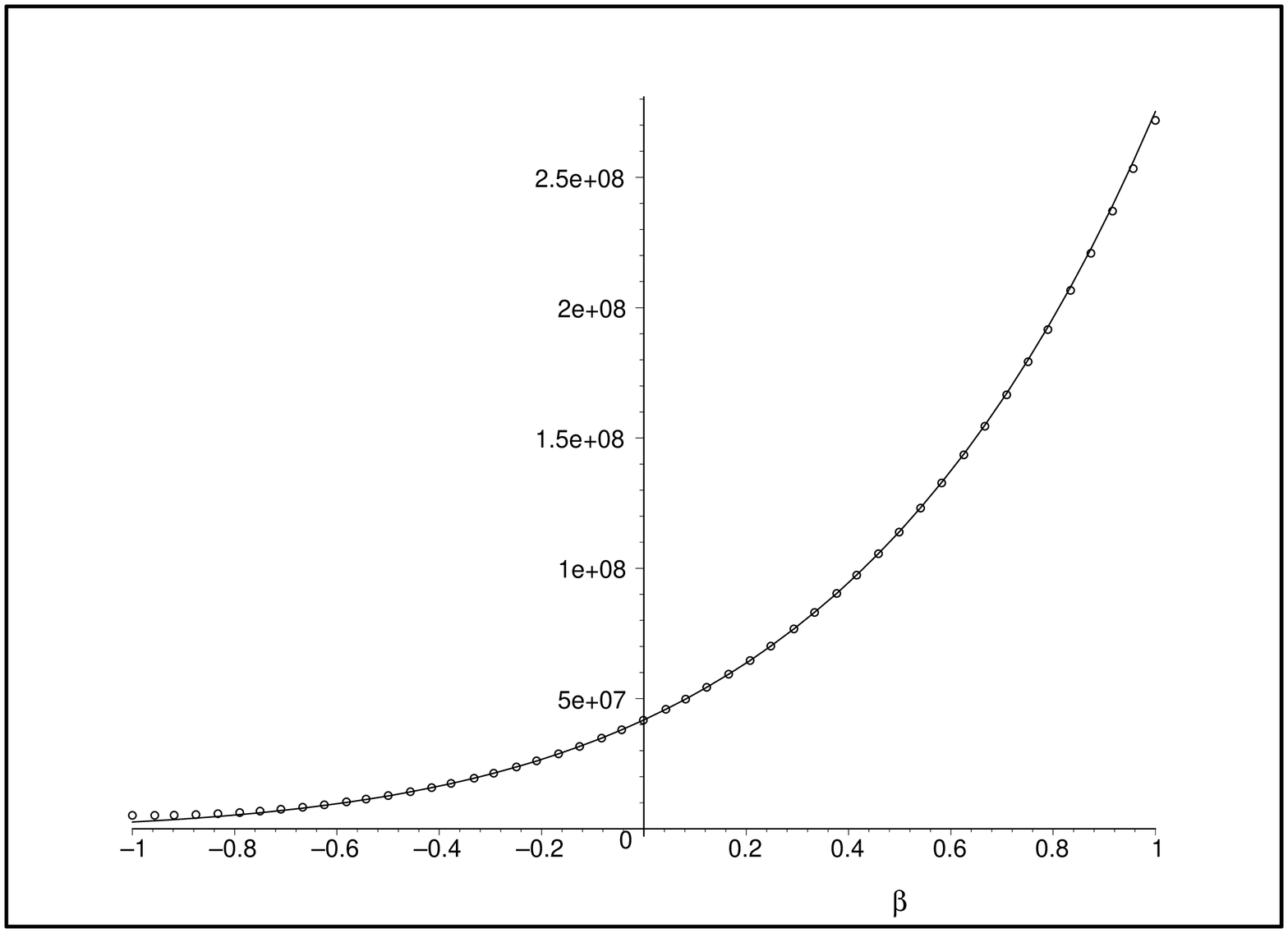}}}
\end{center}
\caption{A comparison of $K_{n}(x)$ (solid curve) and $K^{(8)}(\beta,z)$ (ooo)
for $n=10$ with $\varepsilon=0.01$ and $q=0.34894783$.}%
\label{RegionVIII}%
\end{figure}

\item Region IX: \ $y\approx Y^{-}(z),$ $p<z<1.$%
\begin{align*}
K_{n}(x)  &  \sim K^{(9)}(\beta,z)=\varepsilon^{1/3}\exp\left[  \varepsilon
^{-1}\psi_{0}(z)+\ln\left(  \frac{U_{0}+p}{U_{0}-q}\right)  \beta
\varepsilon^{-1/3}\right] \\
&  \times\frac{1}{2}\frac{\vartheta^{-1/3}}{\sqrt{zU_{0}}}\left[  \lambda
^{+}(\beta,z)\mathrm{Ai}\left(  \vartheta^{2/3}\beta\right)  +\mathrm{i}%
\lambda^{-}(\beta,z)\mathrm{Bi}\left(  \vartheta^{2/3}\beta\right)  \right]  ,
\end{align*}
where $\vartheta(z)=-\Theta(z),$%
\[
\lambda^{\pm}(\beta,z)=\exp\left\{  \frac{2\pi\mathrm{i}\left[  Y^{-}\left(
z\right)  -\beta\varepsilon^{2/3}\right]  }{\varepsilon}\right\}  \pm1.
\]
and $\mathrm{Ai}\left(  \cdot\right)  ,\mathrm{Bi}\left(  \cdot\right)  $ are
the Airy functions (see Figure \ref{RegionIX}).

\begin{figure}[ptb]
\begin{center}
\rotatebox{270} {\resizebox{12cm}{!}{\includegraphics{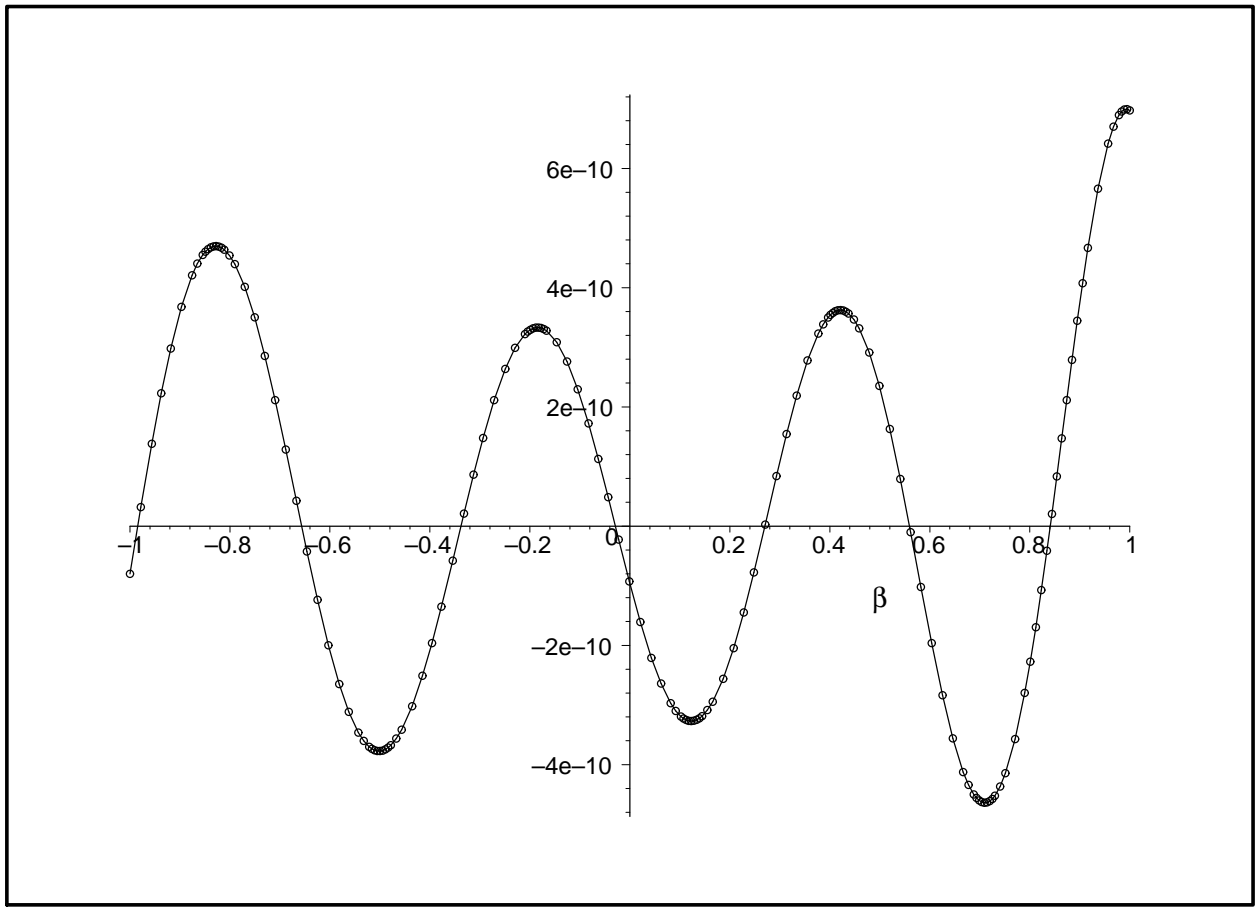}}}
\end{center}
\caption{A comparison of $K_{n}(x)$ (solid curve) and $K^{(9)}(\beta,z)$ (ooo)
for $n=40$ with $\varepsilon=0.02$ and $q=0.74894783$.}%
\label{RegionIX}%
\end{figure}

\item Region X: \ $Y^{-}(z)<y<Y^{+}(z),$ $0<z<1.$%
\[
K_{n}(x)\sim K^{(10)}(y,z)=K^{+}(y,z)+K^{-}(y,z)
\]
(see Figure \ref{RegionX}).

\begin{figure}[ptb]
\begin{center}
\rotatebox{270} {\resizebox{12cm}{!}{\includegraphics{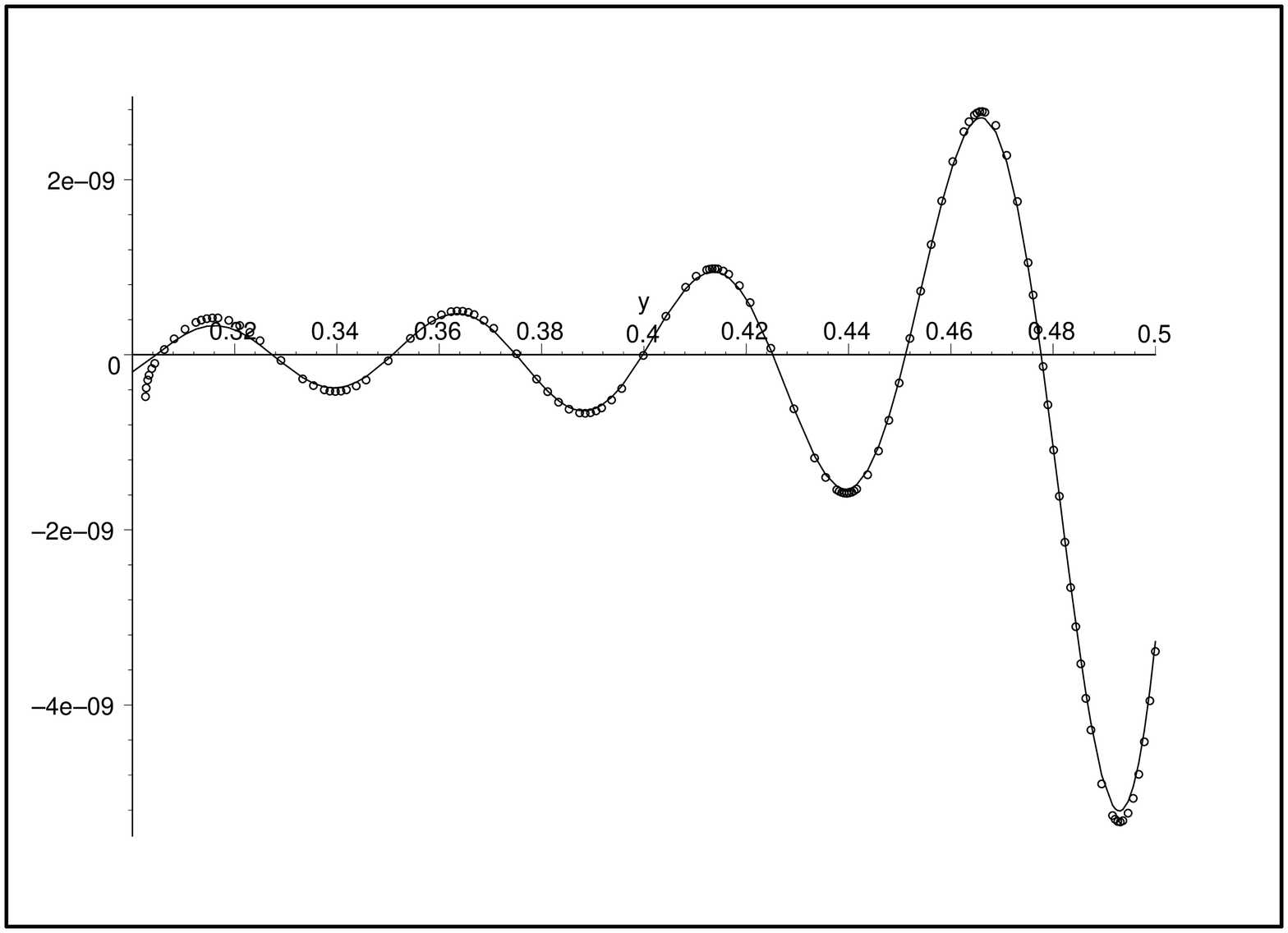}}}
\end{center}
\caption{A comparison of $K_{n}(x)$ (solid curve) and $K^{(10)}(y,z)$ (ooo)
for $n=40$ with $\varepsilon=0.02$ and $q=0.74894783$.}%
\label{RegionX}%
\end{figure}

\item Region XI: \ $n\approx N,$ $n=N-j,$ $j\in\mathbb{Z},$ $0<y<1,$
$y\not \approx q.$%
\[
K_{n}(x)\sim K_{j}^{(11)}(y)=\binom{N}{j}\left(  -p\right)  ^{N-j}\left(
\frac{q}{p}\right)  ^{Ny}\cos\left(  N\pi y\right)  \left(  1-\frac{y}%
{q}\right)  ^{j}+\binom{Ny}{N-j}\left(  \frac{1-y}{q-y}\right)  ^{j+1}%
\]
(see Figure \ref{RegionXI}).

\begin{figure}[ptb]
\begin{center}
\rotatebox{270} {\resizebox{12cm}{!}{\includegraphics{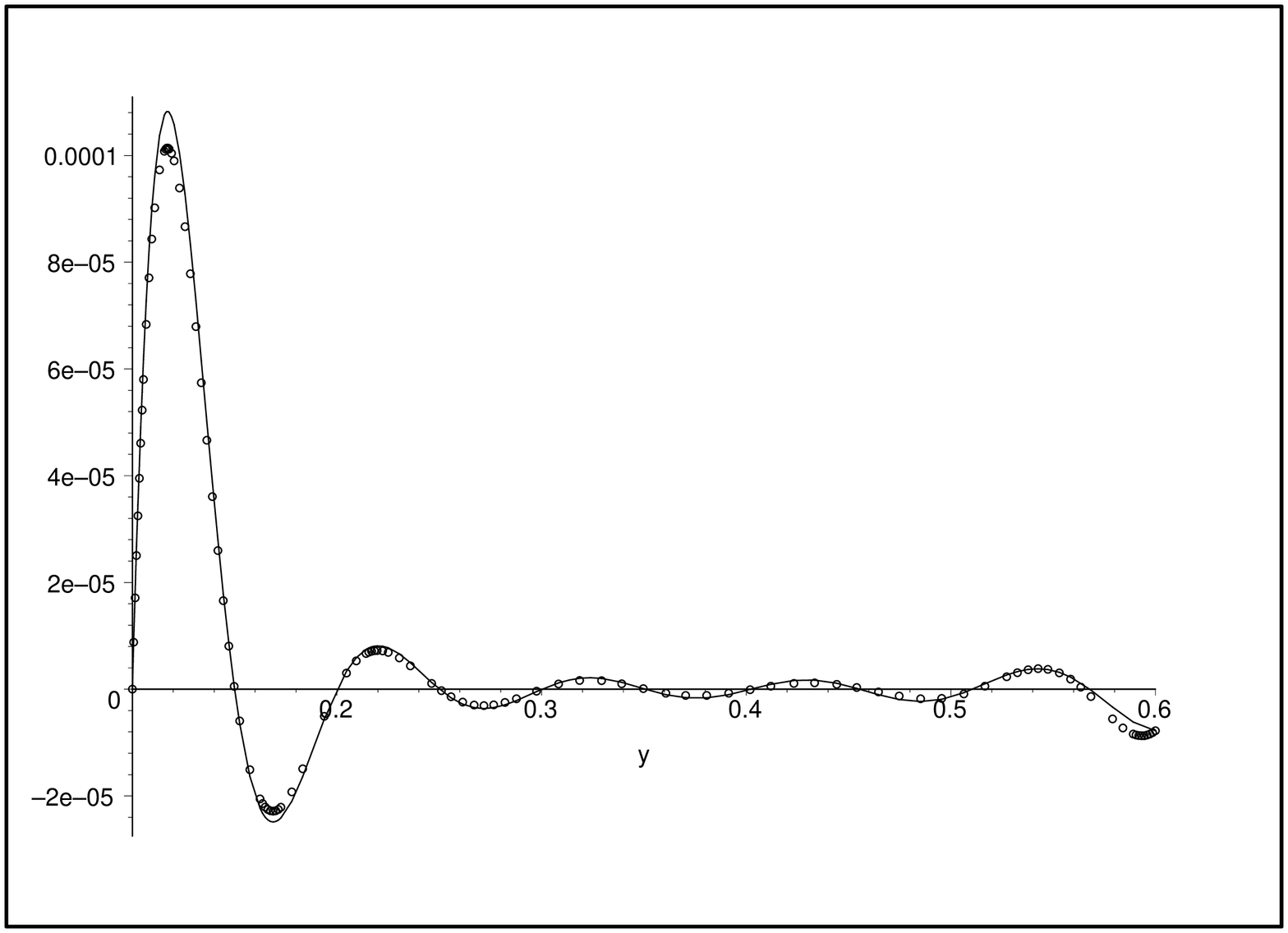}}}
\end{center}
\caption{A comparison of $K_{n}(x)$ (solid curve) and $K_{j}^{(11)}(y)$ (ooo)
for $n=19$ with $\varepsilon=0.05$ and $q=0.74894783$.}%
\label{RegionXI}%
\end{figure}

\item Region XII: \ $n\approx N,$ $y\approx q,$ $y=q+\xi\sqrt{2pq\varepsilon
},$ $\xi=O(1).$%
\begin{gather*}
K_{n}(x)\sim K_{j}^{\left(  12\right)  }(\xi)=\exp\left[  \frac{p\ln
(p)+q\ln(q)}{\varepsilon}+\xi\sqrt{\frac{2pq}{\varepsilon}}\ln\left(  \frac
{q}{p}\right)  -\frac{j}{2}\ln(pq\varepsilon)\right]  \exp\left(  \frac
{\xi^{2}}{2}\right) \\
\times\left[  \frac{1}{j!}\mathrm{D}_{j}\left(  \sqrt{2}\xi\right)
\cos\left(  \frac{p\pi}{\varepsilon}-\xi\pi\sqrt{\frac{2pq}{\varepsilon}%
}\right)  -\frac{1}{\sqrt{2\pi}}\Lambda_{j}\left(  \xi\right)  \sin\left(
\frac{p\pi}{\varepsilon}-\xi\pi\sqrt{\frac{2pq}{\varepsilon}}\right)  \right]
,
\end{gather*}
where the function $\Lambda_{j}:\mathbb{R}\rightarrow\mathbb{R}$ is defined
by
\[
\Lambda_{j}\left(  \xi\right)  =\mathrm{i}^{j+1}\left[  \mathrm{D}%
_{-j-1}\left(  \sqrt{2}\mathrm{i}\xi\right)  +\left(  -1\right)
^{j+1}\mathrm{D}_{-j-1}\left(  -\sqrt{2}\mathrm{i}\xi\right)  \right]
\]
(see Figure \ref{RegionXII}).

\begin{figure}[ptb]
\begin{center}
\rotatebox{270} {\resizebox{12cm}{!}{\includegraphics{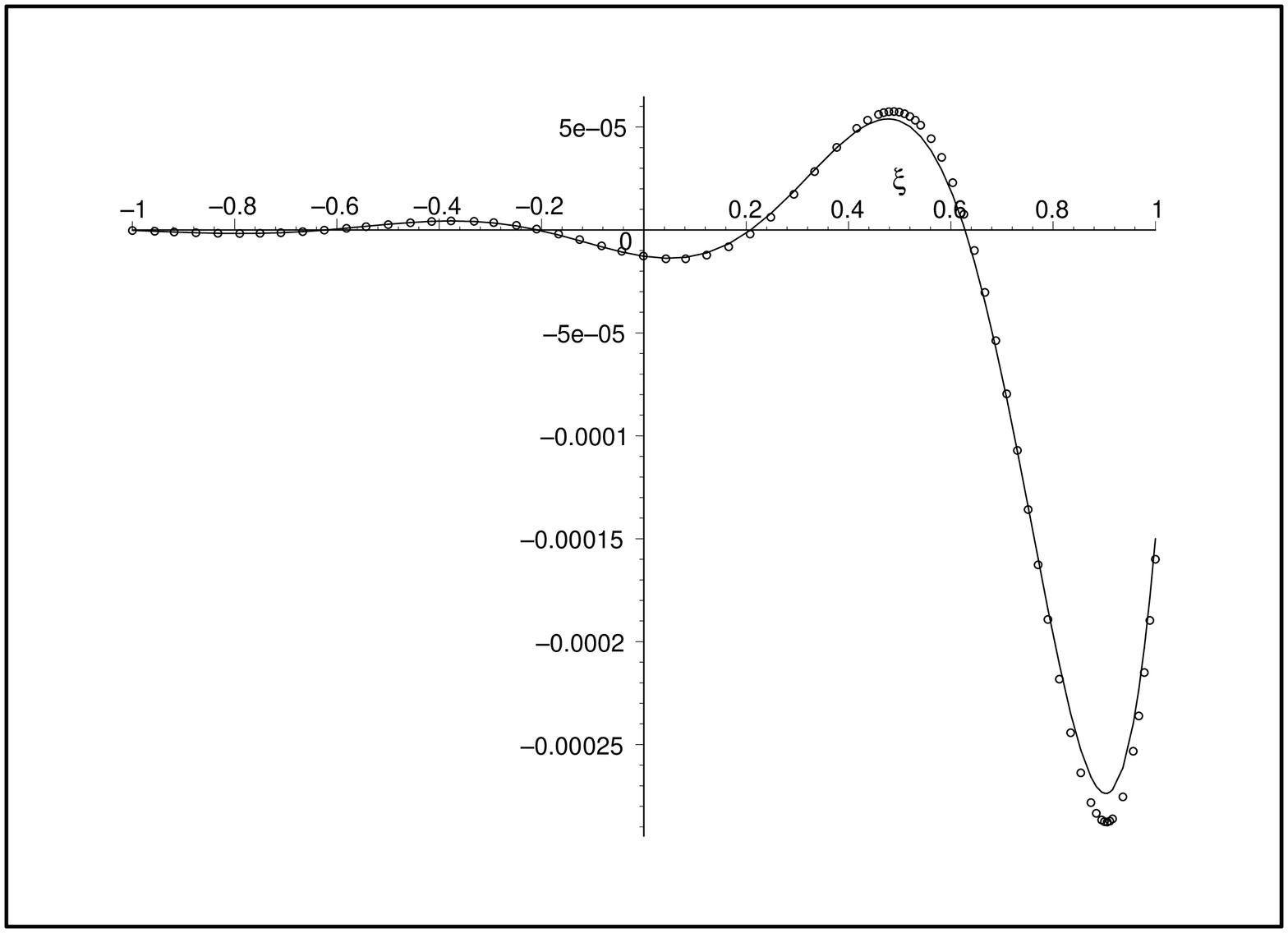}}}
\end{center}
\caption{A comparison of $K_{n}(x)$ (solid curve) and $K_{j}^{\left(
12\right)  }(\xi)$ (ooo) for $n=20$ with $\varepsilon=0.05$ and $q=0.74894783$%
.}%
\label{RegionXII}%
\end{figure}
\end{enumerate}

\end{document}